\documentclass[12pt,reqno]{amsart}
\usepackage{hyperref}

\makeatletter
\def\@tocline#1#2#3#4#5#6#7{\relax
  \ifnum #1>\c@tocdepth 
  \else
    \par \addpenalty\@secpenalty\addvspace{#2}%
    \begingroup \hyphenpenalty\@M
    \@ifempty{#4}{%
      \@tempdima\csname r@tocindent\number#1\endcsname\relax
    }{%
      \@tempdima#4\relax
    }%
    \parindent\z@ \leftskip#3\relax \advance\leftskip\@tempdima\relax
    \rightskip\@pnumwidth plus4em \parfillskip-\@pnumwidth
    #5\leavevmode\hskip-\@tempdima
      \ifcase #1
       \or\or \hskip 1em \or \hskip 2em \else \hskip 3em \fi%
      #6\nobreak\relax
    \dotfill\hbox to\@pnumwidth{\@tocpagenum{#7}}\par
    \nobreak
    \endgroup
  \fi}
\makeatother

\usepackage{tikz,amsthm,amsmath,amstext,amssymb,amscd,epsfig,euscript, mathrsfs, dsfont,pspicture,multicol,graphpap,graphics,graphicx,times,enumerate,subfig,sidecap,wrapfig,color}

\usepackage{ hyperref}
\usepackage{enumerate}

 \numberwithin{equation}{section}

\def\bC{{\mathbb{C}}}

\def\bR{{\mathbb{R}}}
\def\bS{{\mathbb{S}}}

\def\bZ{{\mathbb{Z}}}

\def\bN{{\mathbb{N}}}

\def\cB{{\mathscr{B}}}
\def\cC{{\mathscr{C}}}
\def\cD{{\mathscr{D}}}

\def\cH{{\mathscr{H}}}

\def\cM{{\mathscr{M}}}

\def\cP{{\mathscr{P}}}

\def\cW{{\mathscr{W}}}

\def\one{\mathds{1}}

\def\ve{\varepsilon}

\renewcommand{\d}{{\partial}}
\def\lec{\lesssim}
\def\gec{\gtrsim}

\DeclareMathOperator{\diam}{diam}
\def\dim{\mathop\mathrm{dim}} 					
\def\dist{\mathop\mathrm{dist}} 						
\def\supp{\mathop\mathrm{supp}}					

\newcommand{\ps}[1]{\left( #1 \right)}

\newcommand{\ck}[1]{\left\{#1 \right\}}

\newcommand{\cnj}[1]{\overline{#1}}

\def\XXint#1#2#3{{\setbox0=\hbox{$#1{#2#3}{\int}$ }
\vcenter{\hbox{$#2#3$ }}\kern-.58\wd0}}

\theoremstyle{plain}
\newtheorem{theorem}{Theorem}

\newtheorem{corollary}[theorem]{Corollary}

\newtheorem{lemma}[theorem]{Lemma}

\theoremstyle{definition}

\newtheorem{definition}[theorem]{Definition}
\newtheorem{remark}[theorem]{Remark}

\numberwithin{equation}{section}
\numberwithin{theorem}{section}

\newcommand\eqn[1]{\eqref{e:#1}}
\def\Claim{ {\bf Claim: }}
\newcommand\Theorem[1]{Theorem \ref{t:#1}}
\newcommand\Lemma[1]{Lemma \ref{l:#1}}
\newcommand\Corollary[1]{Corollary \ref{c:#1}}
\newcommand\Definition[1]{Definition \ref{d:#1}}

\newcommand\Remark[1]{Remark \ref{r:#1}}

  \DeclareFontFamily{U}{mathb}{\hyphenchar\font45} 
\DeclareFontShape{U}{mathb}{m}{n}{
      <5> <6> <7> <8> <9> <10> gen * mathb
      <10.95> mathb10 <12> <14.4> <17.28> <20.74> <24.88> mathb12
      }{}
\DeclareSymbolFont{mathb}{U}{mathb}{m}{n}

\DeclareMathSymbol{\toitself}      {3}{mathb}{"FD}  

\begin{document}

\title{Sets of absolute continuity for harmonic measure in NTA domains}

\author{Jonas Azzam}
\address{Departament de Matem\`atiques\\ Universitat Aut\`onoma de Barcelona \\ Edifici C Facultat de Ci\`encies\\
08193 Bellaterra (Barcelona) }
\email{jazzam "at" mat.uab.cat}
\thanks{The author was supported by grants ERC grant 320501 of the European Research Council (FP7/2007-2013) and NSF RTG grant 0838
212.}

\maketitle

\begin{abstract}
We show that if $\Omega$ is an NTA domain with harmonic measure $\omega$ and $E\subseteq \partial\Omega$ is contained in an Ahlfors regular set, then $\omega|_{E}\ll \mathscr{H}^{d}|_{E}$. Moreover, this holds quantitatively in the sense that for all $\tau>0$ $\omega$ obeys an $A_{\infty}$-type condition with respect to $\mathscr{H}^{d}|_{E'}$, where $E'\subseteq E$ is so that $\omega(E\backslash E')<\tau \omega(E)$, even though $\partial\Omega$ may not even be locally $\mathscr{H}^{d}$-finite. We also show that, for uniform domains with uniform complements, if $E\subseteq\partial\Omega$ is the Lipschitz image of a subset of $\mathbb{R}^{d}$, then there is $E'\subseteq E$ with $\cH^{d}(E\backslash E')<\tau \mathscr{H}^{d}(E)$ upon which a similar $A_{\infty}$-type condition holds.
\end{abstract}
\tableofcontents

\section{Introduction}

\subsection{Background}

Given a domain $\Omega\subseteq \bR^{d+1}$ and $E\subseteq \d\Omega$, when do we have $\omega\ll \cH^{d}$ on $E$? In \cite{Ok80}, {\O}ksendal showed that harmonic measure on a simply connected planar domain $\Omega$ is absolutely continuous with respect to to $\cH^{1}$ on $E$ if it is contained in a line $L$. In \cite{KWu82}, Kaufmann and Wu generalized this by showing $L$ could replaced with a bi-Lipschitz curve, and Bishop and Jones in \cite{BJ90} showed absolute continuity occurred inside {\it any} Lipschitz curve. In dimensions larger than two, however, the obvious generalizations of these results are false: In \cite{Wu86}, Wu gives an example of a domain in $\bR^{3}$ that gives positive harmonic measure to a set of Hausdorff dimension 1 in $\bR^{2}$. In spite of this, she proves an analogue of \cite{KWu82} under some mild geometric assumptions. The first involves the notion of uniformity.

\begin{definition}
We say that $\Omega$ is a {\it $C$-uniform domain} if, for every $x,y\in \cnj{\Omega}$ there is a path $\gamma\subseteq\Omega$ connecting $x$ and $y$ such that
\begin{enumerate}
\item the length of $\gamma$ is at most $C|x-y|$ and
\item for $t\in \gamma$, $\dist(t,\d\Omega)\geq \dist(t,\{x,y\})/C$. 
\end{enumerate}
\label{d:uniform}
\end{definition}

Roughly speaking, this says that the domain $\Omega$ has no bottlenecks. The second condition is the following.

\begin{definition}
We say that $\Omega$ satisfies the {\it $C$-interior corkscrew condition} if for all $\xi\in \d \Omega$ and $r\in(0, \diam\d\Omega)$ there is a ball $B(x,r/C)\subseteq \Omega\cap B(\xi,r)$. We say $\Omega$ satisfies the {\it $C$-exterior corkscrew condition} if there is a ball $B(y,r/C)\subseteq B(\xi,r)\backslash\Omega$ for all $\xi\in \d\Omega$ and $r\in(0,\diam \d\Omega)$. 
\label{d:cork}
\end{definition}

It is not hard to show that a $C$-uniform domain satisfies the interior corkscrew condition with constant depending on $C$.

 We can now state the result from \cite{Wu86}. For a domain $\Omega$, we will let $\omega_{\Omega}^{z}$ denote harmonic measure evaluated at a point $z\in\Omega$.

\begin{theorem}\cite{Wu86}
Let $\Omega\subseteq \bR^{d+1}$ be any domain satisfying the exterior corkscrew condition, and let $\Gamma$ be a topological sphere such that $\Gamma^{c}=\Omega_{1}\cup \Omega_{2}$ where $\Omega_{1}$ and $\Omega_{2}$ are disjoint uniform domains for which $\omega_{\Omega_{i}}^{z_{i}}\ll \cH^{d}|_{\Gamma}$ for $i=1,2$ and $z_{i}\in\Omega_{i}$. Then $\omega_{\Omega}^{z}|_{\Gamma\cap \d\Omega}\ll \cH^{d}|_{\Gamma\cap\d\Omega}$ for $z\in\Omega$. 
\label{t:wu}
\end{theorem}

Admissible surfaces $\Gamma$ include, for example, bi-Lipschitz images of $\bS^{d}$ (see \cite[Theorem 10.1]{JerisonKenig}.\\

Under more stringent conditions on the geometry of $\Omega$, one can glean more quantitative information about absolute continuity. The first condition is just the combination of the previous two conditions we've seen so far.

\begin{definition}
A {\it $C$-nontangentially accessible (or $C$-NTA) domain}\footnote{This is not the usual definition of NTA domains, but it is quantitatively equivalent. For example, see \cite{AHMNT}.} $\Omega$ is a $C$-uniform domain satisfying the $C$-exterior corkscrew condition.
\end{definition}

These domains were introduced in \cite{JerisonKenig} by Jerison and Kenig, and they have just enough geometry to guarantee harmonic measure enjoys some useful properties (see \Theorem{JK} below).  

The next assumption is that $\d\Omega$ is Ahlfors regular.

\begin{definition}
A metric space $Z$ is {\it $A$-Ahlfors $d$-regular} if there is $A\geq 1$ so that
\begin{equation}
r^{d}/A\leq  \cH_{Z}^{d}(B_{Z}(x,r))\leq Ar^{d}\mbox{ for all }x\in Z, 0<r<\diam Z
\label{e:regular}
\end{equation}
where $B_{Z}(x,r)$ and $\cH^{d}_{Z}$ denote the open ball in $Z$ of radius $r$ centered at $x$ and the Hausdorff measure on $Z$ respectively.
\end{definition}

In \cite{DJ90}, David and Jerison showed that, under these assumptions, not only are $\omega$ and $\cH^{d}$ mutually absolutely continuous, but quantitatively so, which we make precise in the next theorem.

\begin{theorem}
For all $A,C>1$, integers $d\geq 2$, and $\ve>0$, there are constants $C_{DJ}=C_{DJ}(A,C,d)>0$ and $\delta=\delta(\ve,A,C,d)>0$ such that the following holds. Let $\Omega\subseteq \bR^{d+1}$ be a $C$-NTA domain with an $A$-Ahlfors $d$-regular boundary. Let $\xi_{0}\in \d\Omega$, $r_{0}\in (0,\diam\d\Omega)$, $z_{0}\in \Omega \backslash B(\xi_{0},C_{DJ}r_{0})$, and set  $\omega=\omega_{\Omega}^{z_{0}}$. Then $\omega$ is $A_{\infty}$-equivalent to $\cH^{d}$ on $B_{0}\cap \d\Omega$, meaning whenever $F\subseteq B(\xi,r)\cap E$ with $\xi\in\d\Omega $ and $B(\xi,r)\subseteq B(\xi_{0},r_{0})$, we have\footnote{In \cite{DJ90}, they show that there exists a certain pair $\ve,\delta>0$ so that these conditions hold, but since $\d\Omega$ is Ahlfors regular, this stronger statement can be shown to hold by repeating the arguments in Chapter 5 of \cite{Big-Stein}.}
\begin{enumerate}[(a)]
\item $\omega(F)/\omega(B(\xi,r))<\delta$ implies $\cH^{d}|_{\d\Omega}(F)/ \cH^{d}|_{\d\Omega}(B(\xi,r))<\ve$ and
\item $\cH^{d}|_{\d\Omega}(F)/ \cH^{d}|_{\d\Omega}(B(\xi,r))<\delta$ implies $\omega(F)/\omega(B(\xi,r))<\ve$.
\end{enumerate}
In particular, $\omega\ll \cH^{d}\ll \omega$ on $\d\Omega$. 
\label{t:DJ}
\end{theorem}

This is an improvement over a result of Dahlberg who originally proved this for domains whose boundaries were locally Lipschitz graphs \cite{Dahl77}. In \cite{Badger12}, Badger showed one still has $\cH^{d}\ll \omega$ if instead of \eqn{regular} one only assumes $\cH^{d}|_{\d\Omega}$ is locally finite. He also conjectured that one should still have $\omega\ll \cH^{d}$ in this scenario, but this is false by an example of the author, Mourgoglou, and Tolsa  \cite{AMT15} (which is a refinement of an example of Wolff \cite{Wolff91} that we will discuss below). In fact, locally, the domain $\Omega\subseteq \bR^{d+1}$ constructed in \cite{AMT15} satisfies, for some $A>1>\ve>0$, $A^{-1}r^{d}<\cH^{d}(B(\xi,r)\cap \d\Omega)< Ar^{d-\ve}$ for $\xi\in \d\Omega$ and $r>0$ small; that is, \eqn{regular} just barely fails for $\d\Omega$. 

Finally, we mention some very recent results in the opposite direction, that is, results describing {\it necessary} conditions for absolute continuity. The domains considered by Dahlberg, David and Jerison, and Badger mentioned above all have boundaries that are {\it $d$-rectifiable}, meaning they may be exhausted up to a set of $d$-measure zero by $d$-dimensional Lipschitz graphs, and this is crucial in establishing absolute continuity. The first author, Hofmann, Martell, Mayboroda, Mourgoglou, Tolsa, and Volberg have shown that $\omega|_{E}\ll \cH^{d}|_{E}$ in fact {\it implies} $E$ is a $d$-rectifiable set plus a set of $\omega$-measure zero, and this holds for {\it any} domain $\Omega\subseteq \bR^{d+1}$ for any $d\geq 1$ \cite{AHMMMTV15} (see also \cite{MT15} for a quantitative version of this result). In particular, the result of Bishop and Jones in the plane can be improved: for $\Omega\subseteq \bC$ simply connected, $\omega|_{E}\ll \cH^{1}|_{E}$ for some $E\subseteq \d\Omega$ if and only if $E$ may be covered up to $\omega$-measure zero by Lipschitz curves. 

Hofmann, Le, Martell, and Nystr\"om have shown that if the Poisson kernel for harmonic measure of a uniform domain with $d$-regular boundary satisfies a type of weak-reverse-H\"older inequality, then this implies the boundary is {\it uniformly} rectifiable, and they even prove versions of this for $p$-harmonic measures (see \cite{HLMN15}). See also \cite{ABHM15} and \cite{Mo15}.

\subsection{Main results}
Our results will require the notion of $A_{\infty}$-equivalence on arbitrary sets that may not be Ahlfors regular.

\begin{definition}\label{d:afin}
For a Borel measure $\mu$ in $\bR^{d+1}$ and $E\subseteq \bR^{d+1}$, we will say that {\it $\mu$ is $A_{\infty}$-equivalent to $\cH^{d}$ on $E$} if, for all $\ve>0$, there is $\delta>0$ so that, whenever $F\subseteq B(\xi,r)\cap E$ is a Borel set with $\xi\in E$ and $r>0$,
\begin{enumerate}[(a)]
\item $\mu(F)/\mu(B(\xi,r))<\delta$ implies $\cH^{d}(F)/r^{d}<\ve$ and
\item $\cH^{d}(F)/r^{d}<\delta$ implies $\mu(F)/\mu(B(\xi,r))<\ve$.
\end{enumerate}
We'll say that $\mu$ is $A_{\infty}$-equivalent to $\cH^{d}$ {\it with data depending on} $t_{1},...,t_{n}$ if $\delta$ depends on these as well as $\ve$.
\end{definition}

Observe that if $E=\bR^{d}$, this gives the usual definition of $A_{\infty}$-equivalence.\\

Our first main result generalizes the works of David, Jerison, and Wu mentioned above for the case of NTA domains. Firstly, we remove the requirement that the portion of the boundary in question need be contained in a topological surface as in Wu's theorem. Secondly, we prove an $A_{\infty}$ condition similar to David and Jerison's theorem but on a subset of the boundary contained in a locally $d$-regular set, rather than assuming that the whole boundary is $d$-regular. 

{ 
\begin{theorem}
Let $\Omega\subseteq\bR^{d+1}$ be a $C$-NTA domain. Let $r_{0}\in (0,\diam\d\Omega)$, $\xi_{0}\in\d\Omega$, and $E\subseteq  \d \Omega\cap B(\xi_{0},r_{0})$ be Borel with $\omega(E)/\omega(B(\xi_{0},r_{0}))\geq \rho>0$, where $\omega=\omega_{\Omega}^{z_{0}}$ and $B(z_{0},r_{0}/C)\subseteq   \Omega $. Also suppose there is an $L$-bi-Lipschitz injection $g:E\rightarrow Z$ where  $Z$ is a metric space such that 
\begin{equation}
 r^{d}/A\leq \cH^{d}(B_{Z}(x,r))\leq Ar^{d} \mbox{ for all }x\in g(E), \;\; r\in (0,r_{0}).
 \label{e:almost-ahlfors}
 \end{equation}
Here, $B_{Z}$ is the metric ball in $Z$. Then for all $\tau>0$ there is $E'\subseteq E$ compact and $C^{\pm}$-NTA domains $\Omega_{E'}^{\pm}$ with $C^{\pm}=C^{\pm}(C,d)>0$, such that 
\begin{enumerate}
\item $\omega(E\backslash E')\leq \tau \omega(E)$,
\item $\Omega_{E'}^{-}\subseteq \Omega\subseteq \Omega_{E'}^{+}$,
\item $\Omega_{E'}^{-}\subseteq B(\xi_{0},C^{-}r_{0})$ and $\diam\d\Omega_{E'}^{\pm}\geq r_{0}/C^{-}$,
\item $\d \Omega_{E'}^{\pm}\cap \d \Omega=E'$,
\item $\d \Omega_{E'}^{\pm}$ are $A^{\pm}$-Ahlfors regular with $A^{\pm}$ depending on $A,C,d,L,\rho$ and $\tau$,
\item $\omega$ is $A_{\infty}$-equivalent to $\cH^{d}$ on $E'$ with data depending on $A,C,d,L,\rho$ and $\tau$; in particular, $\omega|_{E'}\ll \cH^{d}|_{E'} \ll \omega|_{E'}$ and $\omega|_{E}\ll \cH^{d}|_{E}$,
\item there is $\delta_{0}>0$ depending on $A,C,d,L,\rho$ and $\tau$ so that $\cH^{d}(E)\geq \cH^{d}(E')\geq \delta_{0}r_{0}^{d}>0$.
\end{enumerate}
\label{t:main}
\end{theorem}
}

\begin{figure}[h]
\scalebox{.42}{\includegraphics{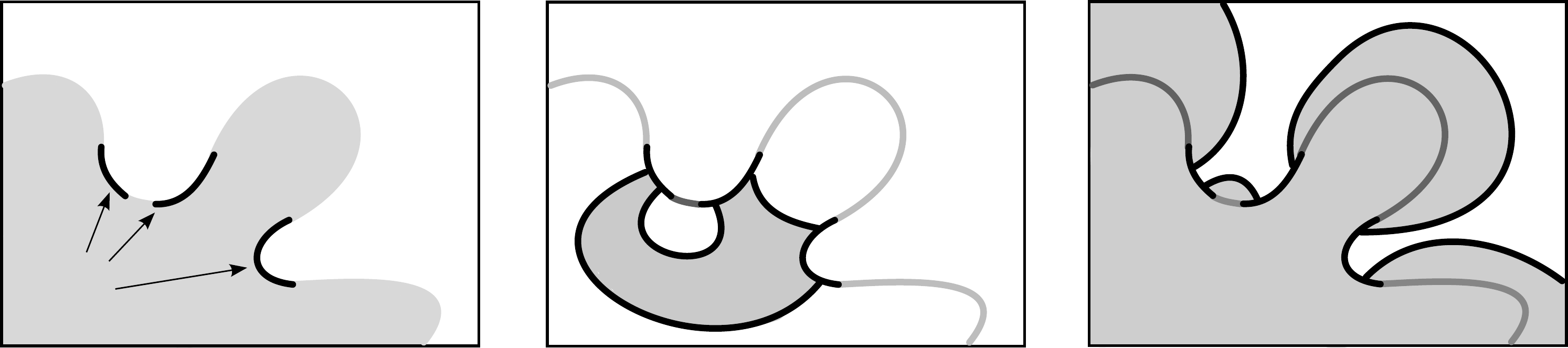}}
\begin{picture}(0,0)(0,0)
\put(-130,0){$\Omega$}
\put(-165,25){$E$}
\put(0,0){$\Omega_{E'}^{-}$}
\put(125,0){$\Omega_{E'}^{+}$}
\end{picture}
\caption{The shaded regions represent $\Omega$, $\Omega_{E'}^{-}\subseteq \Omega$ and $\Omega_{E'}^{+}\supseteq \Omega$. Note that each $\Omega_{E'}^{\pm}$ traces out a portion of the set $E\subseteq \d\Omega$.}\label{f:inandout}
\end{figure}
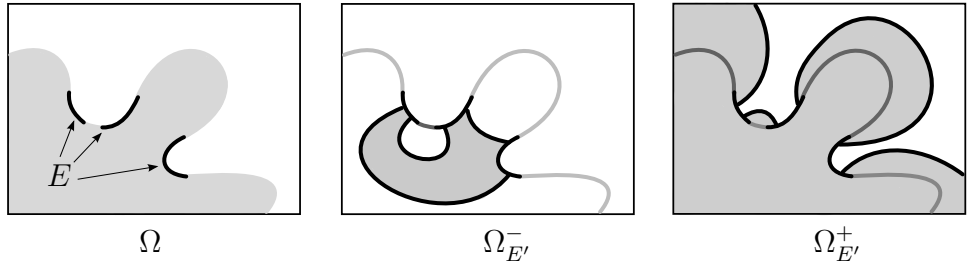

See Figure \ref{f:inandout}. The condition about bi-Lipschitz embedability may seem odd, but one can think instead of $Z$ as being an Ahlfors regular subset of $\bR^{d+1}$ (in which case $g$ is the identity and $L=1$), or of $E$ as a bi-Lipschitz image of a subset of $Z=\bR^{d}$. Observe that we don't assume $\d \Omega$ is Ahlfors regular or even locally $d$-finite as in \Theorem{DJ}. Moreover, we have no a priori assumptions on the set $Z$ other than \eqn{almost-ahlfors}; in \Theorem{wu}, for example, $E$ is assumed to be in a topological surface $\Gamma$ with various conditions, whereas our $Z$ could be totally disconnected. The weak Ahlfors regularity assumption on $Z$ should seem natural in light of the example in \cite{AMT15} mentioned earlier.

Recall that David and Jerison showed NTA domains with Ahlfors regular boundaries have (uniformly) rectifiable boundaries. Hence, the set $E'$ in the lemma is also rectifiable. Moreover, the existence of the NTA domains $\Omega_{E'}^{\pm}$ gives us some quantitative information about the degree of rectifiability of $E$: by specifying $\ve>0$, we can contain all but $\ve$-percent of $E$ inside a uniformly rectifiable set whose parameters are controlled by $\ve$ and the NTA constants of $\Omega$. By an exhaustion argument, we have the following (more digestible) corollary.

\begin{corollary}
Let $\Omega\subseteq \bR^{d+1}$ and $E\subseteq \d\Omega\cap Z$ where $Z\subseteq \bR^{d+1}$ is a set of finite $d$-measure for which $0<\liminf_{r\rightarrow 0} \cH^{d}(B(\xi,r)\cap Z)/r^{d}\leq \limsup_{r\rightarrow 0} \cH^{d}(B(\xi,r)\cap Z)/r^{d}<\infty$ for every $\xi\in E$ (for example, if $Z$ is $d$-regular) and $\omega(E)>0$. Then $E$ contains a rectifiable subset $E'$ of positive $d$-measure such that $\omega(E\backslash E')=0$ and $\omega|_{E'}\ll \cH^{d}|_{E'} \ll \omega|_{E'}$, so in particular, $\omega|_{E}\ll \cH^{d}|_{E}$.
\end{corollary}

The proof of \Theorem{main} relies crucially on a simple lemma about the porosity of sets that have positive doubling measure, which may be of independent interest (for example, see \cite{AM15} for a recent application). The statement requires the definition of dyadic cubes on metric spaces and some extra notation, so rather than stating it here, we strongly encourage the reader to glance at the statement below in \Corollary{porous} and the results in Section \ref{s:doubling} in general.\\

Is there a scenario, or a theorem like \Theorem{main}, where the same results hold instead with the roles of $\cH^{d}$ and $\omega$ reversed? That is, if $\cH^{d}(E)>0$, is there a subset of large $\cH^{d}$-measure upon which $\omega$ and $\cH^{d}$ are $A_{\infty}$-equivalent? Without some further restrictions the answer is a definitive no. If we set
\begin{equation}
\Omega=\bR^{3}_{+}\backslash \bigcup_{n=0}^{\infty}\bigcup_{i,j\in \bZ} B(i2^{-n}e_{1}+j2^{-n}e_{2}+2^{-n}e_{3},2^{-n-10})
\label{e:not2nta}
\end{equation}
then $\Omega$ is NTA but $\omega(\bR^{2})=0$. This is why we can't get $\cH^{d}|_{E}\ll \omega|_{E}$ in \Theorem{main}: we know that we can exhaust almost all of $E$ with respect to {\it harmonic measure}, but a large {\it $d$-measure} portion of $E$ could be hiding somewhere. If we assume $\Omega$ has {\it uniform complement} as well as uniform interior (or {\it doubly uniform}), this in some sense gives $E$ less places to hide and rules out that example. Unfortunately, this is still not enough: In \cite{LVV05}, building off of Wolff's original work in \cite{Wolff91}, Lewis, Verchota, and Vogel construct examples of NTA domains $\Omega\subseteq \bR^{d+1}$, $d\geq 2$, for which there is $F\subseteq \d\Omega$ with $\omega(F^{c})=0$ and $\dim F<d$. Since $\dim \d\Omega\geq d$, by Frostmann's lemma, we can find a set $G\subseteq \d\Omega$ of finite and positive $d$-measure, and $E:=G\backslash F$ is also a set of finite and positive $d$-measure for which $\omega(E)=0$. If we assume $E$ is rectifiable, it turns out this is enough to get an analogue of \Theorem{main}.

{ 
\begin{theorem}
Let $\Omega\subseteq \bR^{d+1}$ be a $C$-uniform domain so that $(\Omega^{c})^{\circ}$ is also $C$-uniform. Let $E\subseteq \d\Omega\cap B(\xi_{0},r_{0})$ where $E$ is the $L$-Lipschitz image of a Borel subset of $[0,r_{0}]^{d}$ such that $\cH_{\infty}^{d}(E)/r_{0}^{d}\geq \rho>0$. Then for all $\eta>0$ there is $E'\subseteq E$ such that, for $B(z_{0},r_{0}/C)\subseteq \Omega$, 
\begin{enumerate}
\item $\cH^{d}_{\infty}(E\backslash E')<\eta \cH^{d}_{\infty}(E)$,
\item $\omega_{\Omega}^{z_{0}}$ is $A_{\infty}$-equivalent to $\cH^{d}$ on $E'$, with constants depending on $C,d,\eta,L$ and $\rho$,
\item  $\omega_{\Omega}^{z_{0}}(E)\geq \delta>0$ for some $\delta$ depending on $C,d,\eta,L$ and $\rho$.
\end{enumerate}
\label{t:uniform-main}
\end{theorem}}
(See the next section for the definition of $\cH^{d}_{\infty}$.) This theorem will follow from the more general \Theorem{in-and-out} below whose statement requires the definition of uniform rectifiability; we will review this in Section \ref{s:in-and-out}. As a corollary, we get a qualitative version of \Theorem{uniform-main}.
{ 
\begin{theorem}
Let $\Omega\subseteq \bR^{d+1}$ be a uniform domain with uniform complement and $\omega=\omega_{\Omega}^{z_{0}}$ be harmonic measure on $\Omega$ for some $z_{0}\in \Omega$. If $E\subseteq \d\Omega$ is a $d$-rectifiable set with $0<\cH^{d}(E)<\infty$, then $\omega\ll \cH^{d}\ll \omega$ on $E\backslash S$ where $\cH^{d}(S)=0$; in particular, $\omega(E)>0$. 
\label{t:uniformtwist}
\end{theorem}}

The example in \cite{AMT15} also happens to be doubly uniform, and thus the set $S$ could very well have positive harmonic measure.

\subsection{Outline}
In Section \ref{s:notation}, we go over some notation and some basic preliminary tools. In Section \ref{s:doubling} we introduce some results about porosity and doubling measures that we will need later on for the special case of harmonic measures. In Section \ref{s:omegaE}, we review and prove some general methods for constructing NTA domains containing a given NTA domain $\Omega$ whose boundaries have prescribed intersections with $\d\Omega$, and under what conditions do they have Ahlfors regular boundaries. 

In Section \ref{s:main} we prove Theorem \ref{t:main}, which follows a similar scheme as the proof of \Theorem{DJ} by David and Jerison. They first showed that, in each ball centered on the boundary of $\Omega$, one can trace out a large portion of $\d\Omega$ by the boundary of a Lipschitz domain (that is, domains whose boundaries are locally $L$-Lipschitz graphs and Alhfors regular boundary), see \cite[Theorem 1]{DJ90}). This is so that they can use Dahlberg's theorem \cite{Dahl77}, which says that $L$-Lipschitz domains  have harmonic measure $A_{\infty}$-equivalent to Hausdorff measure. Knowing this allows them to prove the same property for harmonic measure on $\Omega$ via the maximum principle. In our setting, the domains $\Omega_{E'}^{\pm}$ will play the role of their Lipschitz domains, and we use \Theorem{DJ} instead of Dahlberg's theorem to say harmonic measure is $A_{\infty}$-equivalent to $\cH^{d}$ on these subdomains, after which we repeat the maximum principle argument in \cite{DJ90} (see \Lemma{afineq} below). Hence, the bulk of this section is dedicated to showing how to use the results of Sections \ref{s:doubling} and \ref{s:omegaE} to build the necessary domains $\d\Omega_{E}^{\pm}$.

We don't know whether one can just assume that $\Omega$ satisfies the interior corkscrew condition (recall that some extra topological condition on $\Omega$ is necessary by Wu's example). The NTA assumption is mostly to guarantee that the harmonic measure has some doubling properties (which is used in a critical way) and it helps us construct the Ahlfors regular NTA domains $\Omega_{E'}^{\pm}$ so we can apply \Theorem{DJ} to them. For further discussion on this, see \Remark{generalize}.

In Section \ref{s:in-and-out}, we use the lemmas from Section \ref{s:omegaE} and some results from the theory of uniform rectifiability to prove Theorem \ref{t:uniform-main}. 

\subsection{Acknowledgements}
The author would like to thank John Garnett, Mihalis Mourgoglou, Raanan Schul, and Xavier Tolsa for their very helpful discussions and comments on an early draft, Albert Clop for identifying a mistake, Matthew Badger for pointing out some useful references, and the anonymous referee for his/her critique of the paper.

\section{Notation, Preliminaries, and Harmonic Measure}
\label{s:notation}
We will write $a\lec b$ if there is $C>0$ so that $a\leq Cb$ and $a\lec_{t} b$ if the constant $C$ depends on the parameter $t$. We write $a\sim b$ to mean $a\lec b\lec a$ and define $a\sim_{t}b$ similarly. 

In a metric space $Z$, we will denote the distance between points $x,y\in Z$ as $|x-y|$. For sets $A,B\subseteq Z$, we let 
\[\dist(A,B)=\inf\{|x-y|:x\in A,y\in B\}, \;\; \dist(x,A)=\dist(\{x\},A),\]
and 
\[\diam A=\sup\{|x-y|:x,y\in A\}.\]
Set for a set $A\subseteq Z$, let $w_{d}$ be the volume of the unit ball in $\bR^{d}$ and define
\[\cH^{d}_{\delta}(A)=w_{d}\inf\ck{\sum r_{i}^{d}: A\subseteq \bigcup B(x_{i},r_{i}),x_{i}\in\bR^{d}}.\]
where $\omega_{d}$ is the volume of the unit ball in $\bR^{d}$. We define the {\it $d$-dimensional Hausdorff measure} as
\[\cH^{d}(A)=\lim_{\delta\downarrow 0}\cH^{d}_{\delta}(A)\]
and the {\it $d$-dimensional Hausdorff content} as $\cH^{d}_{\infty}(A)$. See \cite[Chapter 4]{Mattila} for more details. 

We will let $B(x,r)$will not denote the open ball of radius $r$ centered at $x$. In this paper, we will be working in either $\bR^{d}$, $\bR^{d+1}$, or a metric space $Z$, and we won't distinguish our notation for $|x-y|$ or $B(x,r)$ in these cases when it is clear from the context what we mean; otherwise, we will let $B_{Z}(x,r)$ denote the ball in $Z$ and $\cH^{d}_{Z}$ Hausdorff measure on $Z$.  Also define  $\lambda B(x,r)=B(x,\lambda r)$ and $\one_{A}$ to be the function identically one on $A$ and zero elsewhere.\\

For $n\in\bZ$, a {\it $d$-dimensional dyadic cube $Q$ of side length $2^{n}$} in $\bR^{d}$ is a $d$-fold Cartesian product of closed intervals of the form $[i2^{n},(i+1)2^{n}]$, where $i\in\bZ$, and we will denote the side length by $\ell(Q)=2^{n}$. We will write $\lambda Q$ for the cube of the same center as $Q$ and edges parallel to the coordinate axes but side length $\lambda\ell(Q)$.

\begin{definition}[Whitney Cubes]
For an open set $\Omega\subseteq\bR^{d+1}$ and $K>1$, we will denote by $\cW_{K}(\Omega)$ the set of maximal dyadic cubes $Q\subseteq \Omega$ such that $K Q\cap \Omega^{c}=\emptyset$. These cubes have disjoint interiors and can be easily shown to satisfy the following properties:
\begin{enumerate}
\item $\frac{K-1}{2}\ell(Q)\leq \dist(x,\Omega^{c})\leq (1+K)\diam Q$ for all $x\in Q$,
\item $(\frac{K-1}{2}-\sqrt{d+1}\frac{\lambda-1}{2})\ell(Q)\leq \dist(x,\Omega^{c})\leq  (1+K+(\lambda-1)/2)\diam Q$ for all $x\in \lambda Q$ if $\lambda\in [1,K)$ is close enough to $1$ (depending on $d$ and $K$)
\item If $Q,R\in \cW_{K}(\Omega)$ intersect, then $\ell(Q)\sim_{K,d}\ell(R)$.
\item $\sum_{Q\in \cW_{k}(\Omega}\one_{\lambda Q}\lec_{K,d}\one_{\Omega}$ for $\lambda\in(1,K)$
\end{enumerate}
We will just write $\cW_{3}(\Omega)$ as $\cW(\Omega)$.
\label{d:Whitney}
\end{definition}
We will say $Q,R\in \cW(\Omega)$ are {\it adjacent} if $Q\cap R\neq\emptyset$ and write $Q\sim R$. Also, let $P_{Q,R}$ denote the shortest path $Q=Q_{0},...,Q_{n}=R$ of Whitney cubes such that $ Q_{j}\sim Q_{j+1}$ for $j=0,...,n-1$ and define $d_{\Omega}(Q,R)=n+1$. With the definition of Whitney cubes and this notation, we can now state an equivalent characterization of $C$-uniformity that we will need later.

\begin{theorem}[Alternate characterization of uniform domains] A domain $\Omega$ is uniform if and only if it satisfies the interior corkscrew condition and there is $N_{\Omega}:[0,\infty)\toitself$ increasing such that,
\begin{equation}
d_{\Omega}(Q,R)\leq N_{\Omega}(\dist(Q,R)/\min\{\ell(Q),\ell(R)\}) \mbox{ for all }Q,R\in \cW(\Omega).
\label{e:d<N}
\end{equation}
\label{t:AHMNT}
\end{theorem}

\begin{remark}
There are a few papers all giving different yet equivalent definitions of uniform domains. A proof of \Theorem{AHMNT} is given in \cite{AHMNT}; there they instead work with the so-called {\it Harnack chain condition}, which is quantitatively equivalent to the characterization in the above theorem. 
\end{remark}

\begin{remark}
As mentioned in the introduction, a $C$-uniform domain $\Omega$ automatically satisfies the interior corkscrew condition, with constant depending on $C$. For the sake of cleanliness, we will assume that all $C$-NTA domains also satisfy the exterior and interior corkscrew conditions {\it with the same constant} $C$ (which can be arranged by increasing the value $C$ depending only on some universal constant).
\label{r:icork}
\end{remark}

Bounded NTA domains $\Omega$ are regular in the sense of Wiener, so given a continuous $f$ on $\d\Omega$, one can use the Perron method to find $u_{f}$ harmonic, continuous up to the boundary, and equal to $f$ on the boundary as in \cite[Section 2.8]{GilbargTrudinger}. Then, given $z\in \Omega$, one defines harmonic measure via the Riesz representation theorem as the Radon measure $w^{z}_{\Omega}$ so that $\int_{\Omega} fdw^{z}=u_{f}(z)$. For unbounded NTA domains, the situation is more complicated, but given a bounded continuous $f$ on $\d\Omega$ we can still find a bounded harmonic $u_{f}$ continuous up to $\d\Omega$ and equal to $f$ there, and thus we can define harmonic measure the same way; we refer the reader to \cite[Chapter 5]{LesterHelms}, particularly pages 206-7, Theorem 5.4.2, and page 217. 

We recall a few basic results from \cite{JerisonKenig}\footnote{Note that in \cite{JerisonKenig}, they assume their NTA domains are bounded domains, but this is only so that they can guarantee the existence of harmonic measure. Now that we know existence also holds for unbounded domains, the results carry over to this setting with identical proofs.}:

\begin{theorem}[Local properties of harmonic heasure]
Let $\Omega$ be a $C$-NTA domain, $\omega_{\Omega}^{z_{0}}$ be harmonic measure evaluated at $z_{0}\in \Omega$, $\xi\in \d\Omega$, $r\in (0,\diam\d\Omega)$, and let $E\subseteq B(\xi,r)\cap \d\Omega$ be Borel.
\begin{enumerate}

\item \cite[Lemma 4.11]{JerisonKenig}  If $B(z,r/C)\subseteq \Omega \cap B(\xi,r)$ and $z_{0}\in \Omega\backslash B(\xi,2r)$, then
\begin{equation}
w^{z_{0}}_{\Omega}(E)/w^{z_{0}}_{\Omega}(B(\xi,r))\sim_{C,d} w^{z}_{\Omega}(E).
\label{e:wratio}
\end{equation}
\item \cite[Lemma 4.2]{JerisonKenig} If $B(z,r/C)\subseteq \Omega\cap B(\xi,r)$, then
\begin{equation}
w^{z}_{\Omega}(B(\xi,r))\gec_{C,d} 1.
\label{e:wbig}
\end{equation}
\item ({\it Harnack's inequality}) If $x\in Q\in \cW(\Omega)$ and $y\in R\in \cW(\Omega)$ and $\dist(Q,R)/\min\{\ell(Q),\ell(R)\}\leq \Lambda$, then for any Borel set $A\subseteq \d\Omega$,
\begin{equation}
w^{x}_{\Omega}(A)\lec_{C,d,\Lambda}w^{y}_{\Omega}(A).
\label{e:wharnack}
\end{equation}
\item ({\it Local doubling property},  \cite[Lemma 4.9]{JerisonKenig}) If $z_{0}\in \Omega\backslash B(\xi_{0},2r_{0})$ for some $\xi_{0}\in \d\Omega$ or $B(z_{0},r_{0}/C)\subseteq\Omega$, then for all $\xi\in \d\Omega$ and $r>0$ with $B(\xi,2r)\subseteq B(\xi_{0},r_{0})$, we have 
\begin{equation}
\omega_{\Omega}^{z_{0}}(B(\xi,2r))\lec_{C,d,M_{0}} \omega_{\Omega}^{z_{0}}(B(\xi,r)).
\label{e:wdub}
\end{equation}
\end{enumerate}
\label{t:JK}
\end{theorem}

As originally stated in \cite{JerisonKenig}, the constant in \eqn{wdub} is also allowed to depend on $z_{0}$, but an inspection of the proof shows that, so long as $z_{0}\in \Omega\backslash B(\xi_{0},2r_{0})$, this inequality holds independent of $z_{0}$, and this implies the case of $B(z_{0},r_{0}/C)\subseteq \Omega$ by Harnack's inequality.

We now have enough tools to demonstrate the role the approach regions $\Omega_{E}^{\pm}$. Again, this type of argument appears in many sources and is rooted in complex analysis and the study of nontangential limits of harmonic functions, as well as the study of harmonic measure in NTA domains, so the lemma below should be considered review. For a survey of this history, see the introduction to \cite{JerisonKenig}.

\begin{lemma}
Let $\Omega$ be a $C$-NTA domain and $\omega=\omega_{\Omega}^{z_{0}}$ where $B(z_{0},r_{0}/C)\subseteq \Omega$ for some $r_{0}\in (0,\diam\d\Omega)$. Suppose $E\subseteq \d\Omega\cap B(\xi_{0},r_{0})$ for some $\xi_{0}\in \d\Omega$ and there are domains $\Omega_{E}^{\pm}$ so that
\begin{enumerate}
\item $\Omega_{E}^{-}\subseteq \Omega\subseteq \Omega_{E}^{+}$,
\item $\d\Omega_{E}^{\pm}\cap \d\Omega\supseteq E$,
\item $\diam \d\Omega_{E}^{\pm}\geq r_{0}/C^{-}$,
\item $\d\Omega_{E}^{\pm}$ is $A^{\pm}$-Ahlfors regular.
\end{enumerate}
Then $\omega$ is $A_{\infty}$-equvalent to $\cH^{d}$ on $E$.
\label{l:afineq}
\end{lemma}

\begin{proof}
Let $E$, $\Omega_{E}^{\pm}$, $\xi\in E$ and $r>0$, and set $B=B(\xi,r)$. By a covering argument, since $E\subseteq B(\xi_{0},r_{0})$, we can assume without loss of generality that $r<qr_{0}$ where $q>0$ will be determined later. Let $B(z,r_{0}/(C^{-})^{2})\subseteq B(\xi,r_{0}/C^{-})\cap \Omega_{E}^{-}$, which exists by the corkscrew condition for $\Omega_{E}^{-}$. Then $\dist(z,\d\Omega_{E}^{-})\geq r_{0}/(C^{-})^{2}$, so if $q^{-1}>2C_{DJ}(C^{-})^{2}$, then $r<qr_{0}$ implies $z\not\in B(\xi,C_{DJ}r)$. By \Theorem{DJ}, we know that $\omega_{\Omega_{E}^{-}}^{z}$ is $A_{\infty}$-equivalent to $\cH^{d}$ on $B\cap \d\Omega_{E}^{-}$, meaning for every $\ve>0$, there is $\delta_{-}=\delta_{-}(\ve,C^{-},d,A^{-})$ so that if $F\subseteq B\cap E$, then 
\begin{enumerate}[(a)]
\item  $\omega_{\Omega_{E}^{-}}^{z}(F)<\delta_{-} \omega_{\Omega_{E}^{-}}^{z}(B)$ implies $\cH^{d}(F)<\ve \cH^{d}(B\cap \d \Omega_{E}^{-})$ and  
\item $\cH^{d}(F)<\delta \cH^{d}(B\cap \d \Omega_{E}^{-})$ implies $\omega_{\Omega_{E}^{-}}^{z}(F)<\ve \omega_{\Omega_{E}^{-}}^{z}(B)$.
\end{enumerate}
Let $\ve,\ve',\delta>0$, and $F\subseteq E\cap B$ where $B=B(\xi,r)$ with $\xi\in E$ and $r<qr_{0}$ (where $\ve',\delta$ and $q$ will be determined later), and assume $\omega(F)<\delta \omega(B)$. Set $\delta'=\delta_{-}(\ve',C^{-},d,A^{-})$. Pick $B(z',r/C^{-})\subseteq B\cap \Omega_{E}^{-}$. For $q^{-1}<2C$, since $B(z_{0},r_{0}/C)\subseteq \Omega$,
\[
|z_{0}-\xi|\geq  r_{0}/C>2r\]
and we have $|z-\xi|\geq C_{DJ}r>2r$ as well, so that we can apply \Theorem{JK} twice along with the maximum principle to obtain
\[\frac{\omega_{\Omega_{E}^{-}}^{z}(F)}{\omega_{\Omega_{E}^{-}}^{z}(B)}
\sim_{C^{-}} \omega_{\Omega_{E}^{-}}^{z'}(F)
\leq \omega_{\Omega}^{z'}(F)
\sim_{C} \frac{\omega(F)}{\omega(B)}<\delta.\]
Thus, for $\delta$ small enough, we have $\omega_{\Omega_{E}^{-}}^{z}(F)<\delta' \omega_{\Omega_{E}^{-}}^{z}(B)$, which implies $\cH^{d}(F)<\ve' \cH^{d}(B\cap \d \Omega_{E}^{-})\leq A^{-} \ve' r^{d}$, and for $\ve'$ small enough, this  implies $\cH^{d}(F)< \ve r^{d}$. 

Conversely, let $\ve,\ve',\delta>0$ (the latter two will be decided soon) and  $F\subseteq E\cap B$ with $\cH^{d}(F)< \delta r^{d}$, were we will decide $\delta$ in a moment. Let $B(z,r/(C^{+})^{2}) \subseteq B(\xi,r/C^{+})\cap\Omega_{E}^{+}$. Again, $q^{-1}>2C_{DJ}(C^{+})^{-2}$ implies $z\not\in B(\xi,C_{DJ} r)$. Let $\ve'>0$ and pick $\delta'$ so that our $A_{\infty}$ condition on $\omega_{\Omega^{+}}^{z}$ holds for $\ve'$ and $\delta'$ on $B\cap \d\Omega_{E}^{+}$. Then $\cH^{d}(F)<\delta r^{d}\leq A^{+}\delta \cH^{d}(B\cap \d\Omega_{E}^{+})$, and so for $\delta$ small enough, $\cH^{d}(F)<\delta'\cH^{d}(B\cap \d\Omega_{E}^{+})$. Again by the maximum principle, the $A_{\infty}$ condition, and \Theorem{JK},
\[\frac{\omega(F)}{\omega(B)}\sim_{C} \omega_{\Omega}^{z}(F) \leq \omega_{\Omega_{E}^{+}}^{z}(F)<\ve' \omega_{\Omega_{E}^{+}}^{z}(B)<\ve'\]
Picking $\ve'$ small enough guarantees $\omega(F)/\omega(B)<\ve $.

\end{proof}

\section{``Cubes'' and Carleson packing conditions on porosity}
\label{s:doubling}

In this section, we will review and develop some tools that will help us find the desired set $E'$ in \Theorem{main}. The material for this section holds in more generality than just harmonic measure on NTA domains, but for doubling measures on metric measure spaces (if it bugs the reader, s/he can imagine all the measures below are just harmonic measure). We start by introducing the notion of ``dyadic cubes" for a metric space. We'll use the construction of Hyt\"onen and Martikainen from \cite{HytMart12}, which refines the originals of Christ \cite{Christ-T(b)} and David \cite{David88}. We will abuse notation by letting $|x-y|$ denote the metric distance between points $x$ and $y$ and $B(x,r)$ again denote the ball centered at $x$ of radius $r$ in the given space.

\begin{theorem}
For $c_{0}<1/1000$, the following holds. Let $c_{1}=1/500$ and $\Sigma$ be a metric space. For each $n\in\bZ$ there is a collection $\cD_{n}$ of ``cubes,'' which are Borel subsets of $\Sigma$ such that
\begin{enumerate}
\item $\Sigma=\bigcup_{\Delta\in \cD_{n}}\Delta$ for every $n$,
\item if $\Delta,\Delta'\in \cD=\bigcup \cD_{n}$ and $\Delta\cap\Delta'\neq\emptyset$, then $\Delta\subseteq \Delta'$ or $\Delta'\subseteq \Delta$,
\item for $\Delta\in \cD_{n}$, there is $\zeta_{\Delta}\in X_{n}$ so that if $B_{\Delta}=B(\zeta_{\Delta},5c_{0}^{n})$, then
\[c_{1}B_{\Delta}\subseteq \Delta\subseteq B_{\Delta}.\]
\end{enumerate}
\label{t:Christ}
\end{theorem}

For $\Delta\in \cD_{n}$, define $\ell(\Delta)=5c_{0}^{n}$, so that $B_{\Delta}=B(\zeta_{\Delta},\ell(\Delta))$. Note that for $\Delta\in \cD_{n}$ and $\Delta'\in \cD_{m}$, we have $\ell(\Delta)/\ell(\Delta')=c_{0}^{n-m}$.

For $\lambda \leq 1$, define
\[\lambda \Delta = 
\{\xi\in \Delta: \dist(\xi,\Sigma\backslash \Delta)>(1-\lambda)\ell( \Delta)\}. \]

Let $\mu$ be a doubling measure on a metric space $\Sigma$, meaning $\mu(B(\xi,2r))\leq C_{\mu} \mu(B(\xi,r))$ for all $\xi\in \Sigma$ and $r>0$. If $E\subseteq \Sigma$ is a $\delta$-porous set (meaning for every $\xi\in E$ and $r>0$ there is $B(\xi',\delta r)\subseteq (\Sigma\backslash E)\cap B(\xi,r)$), then $\mu(E)=0$. This follows from the fact that the Lebesgue differentiation theorem still holds for doubling measures, and $\delta$-porosity implies $\mu(E\cap B(\xi,r))/\mu(B(\xi,r))\leq 1-C_{\delta,\mu}<1$ for all $\xi\in E$ and $r>0$. Thus, a set of positive measure can't be porous inside every ball centered on $E$. In this section, we will quantify how many cubes $\Delta$ there are inside a given cube $\Delta_{0}$ for which a set $E$ is too porous near $\Delta$, and we will give this control in terms of a so-called {\it Carleson packing condition}. We will then use this condition to trim down the set $E$ to a slightly smaller set $E'$ such that every point in $E'$ is contained in at most a bounded number of cubes that are porous for $E$. 

However, we need to be even more careful for our applications later: $\omega_{\Omega}^{z_{0}}$ is globally doubling {\it with doubling constant depending on $z_{0}$}, and we'd like the constants in our results not to have this dependence. By \eqn{wdub}, however, we can guarantee that $\omega_{\Omega}^{z_{0}}$ is doubling locally with constant independent of $z_{0}$ so long as $z_{0}$ avoids that portion of the boundary. Our next lemma, for example, is well known for the case of doubling measures, but we need to alter it a bit to account for the local doubling case.

\begin{lemma}
Let $\Sigma$ be a metric space and $\cD_{n}$ the ``cubes" constructed in \Theorem{Christ}. Let $c_{0}<c_{1}/4$ and $\mu$ be measure on $\Sigma$ such that, for some $\Delta_{0}\in \cD$, $\mu(B(\xi,2r))\leq C_{\mu}\mu(B(\xi,r))$ for $\xi\in 4B_{\Delta_{0}}$ and $0<r\leq \ell(\Delta_{0})$. There are $t_{0},\alpha>0$ (depending on $C_{\mu}$) such that for $t\in (0,1)$ and $\Delta\subseteq \Delta_{0}$, $\mu(\Delta\backslash (1-t)\Delta) \leq t_{0} t^{\alpha}\mu(\Delta)$ for some $\alpha>0$.
\label{l:local-dub}
\end{lemma}

Note that $C_{\mu}$ is not necessarily the doubling constant of $\mu$, only for those particular values of $\xi$ and $r$ in the lemma. This lemma can be obtained by carefully reading the proof in \cite{Christ-T(b)}, but we will provide a proof for the reader's convenience in the appendix.\\

We now give our first lemma that helps quantify how porous a set can be.

\begin{lemma}
Let $\mu$ be a Borel measure, $\cD$ the cubes for $\Sigma=\supp \mu$ with constant $c_{0}$ and $c_{1}$, and $E\subseteq \Delta_{0}\in \cD$ be Borel. Let $M_{0}>1$ and suppose $\mu$ has the property that, for all $\xi\in 4M_{0}B_{\Delta_{0}}$ and $r\in(0,4M_{0}\ell(\Delta_{0}))$, $\mu(B(\xi,2r))\leq C_{\mu}\mu(B(\xi,r))$ for some $C_{\mu}>1$. Let $M<M_{0}$ and $\cW(E^{c})$ denote the maximal cubes $\Delta\in \cD$ such that $MB_{\Delta}\cap E=\emptyset$. For $\beta>0$, set
\begin{equation}
\lambda_{M,\beta} (\Delta)=\sum_{\Delta' \in \cW(E^{c})\atop \Delta'\subseteq MB_{\Delta}} \ps{\frac{\ell(\Delta')}{\ell(\Delta)}}^{\beta}.
\label{e:lambda}
\end{equation}
We set $\lambda_{M,\beta}(\Delta)=0$ if it is an empty sum. If $\beta>\beta_{0}:=\log_{2} C_{\mu}$, then for all $\Delta_{1}\subseteq \Delta_{0}$, 
\begin{equation}
\sum_{\Delta\subseteq \Delta_{1} \atop \Delta\cap E\neq\emptyset} \lambda_{M,\beta}(\Delta) \mu(\Delta ) \lec_{\mu,M,c_{0},c_{1},\beta} \mu(\Delta_{1}).
\end{equation}
\label{l:porous}
\end{lemma}

To avoid some double subscripts, we write $\lec_{\mu}$ to mean that the implied constant depends on the doubling constant $C_{\mu}$.

\begin{proof}
\Claim When $\Delta'\subseteq MB_{\Delta}$ and $\Delta\subseteq \Delta_{0}$,
\begin{equation}
(\ell(\Delta')/\ell(\Delta))^{\beta_{0}}\lec_{M,C_{\mu}} \mu(\Delta' )/\mu(\Delta).
\label{e:b0}
\end{equation}

Let $\Delta'\subseteq MB_{\Delta}$ and let $N$ be such that 
\begin{equation}
\label{e:2^N}
2^{N}c_{1}\ell(\Delta')>2M\ell(\Delta)\geq 2^{N-1}c_{1}\ell(\Delta').
\end{equation}
Then $2^{N}c_{1}B_{\Delta'}\supseteq MB_{\Delta}$, and $2^{N}<\frac{4M\ell(\Delta)}{c_{1}\ell(\Delta')}$, so that 
\[N<\log_{2}\ps{\frac{4M\ell(\Delta)}{c_{1}\ell(\Delta')}}.\]
Thus
\begin{align*}
\mu(\Delta')
& \geq \mu(c_{1}B_{\Delta'})
\geq C_{\mu}^{-N}\mu(2^{N}c_{1}B_{\Delta'}) 
 \geq C_{\mu}^{-N}\mu(MB_{\Delta})\\
 & \geq C_{\mu}^{\log_{2}\frac{c_{1}}{4M}}\ps{\frac{\ell(\Delta')}{\ell(\Delta)}}^{\log_{2}C_{\mu}}\mu(\Delta).
\end{align*}
This proves the claim. From now on, we write $\lambda=\lambda_{M,\beta}$ with $\beta>\beta_{0}$.

For fixed $\Delta_{1}\subseteq \Delta_{0}$, $\Delta'\subseteq MB_{\Delta_{1}}$ and $\Delta'\in \cW(E^{c})$, set
\[\cM_{n}(\Delta')=\{\Delta\in \cD_{n}: \Delta\subseteq \Delta_{1}, \Delta\cap E\neq\emptyset, \Delta'\subseteq MB_{\Delta}\}.\]
We prove here a few properties of this set:
\begin{enumerate}
\item[$(\dag)$] If $\Delta\in \cM_{n}(\Delta')$, then $\ell(\Delta')\leq \ell(\Delta)$, and in particular, $\cM_{n}(\Delta')\neq\emptyset$ only when $5c_{0}^{n}\geq \ell(\Delta')$. 

To see this, note that since $\Delta'\in \cW(E^{c})$, $\Delta\cap E\neq\emptyset$, and $\Delta'\subseteq MB_{\Delta}$, if $\xi\in \Delta\cap E$, we have
\begin{multline*}M\ell(\Delta') - \ell(\Delta) 
\leq \dist(\zeta_{\Delta'},E)-\ell(\Delta)
\leq |\zeta_{\Delta'}-\xi|-\ell(\Delta)|\\
\leq |\zeta_{\Delta'}-\xi|-|\xi-\zeta_{\Delta}|
\leq |\zeta_{\Delta'}-\zeta_{\Delta}|
\leq M\ell(\Delta)
\end{multline*}
so $\ell(\Delta')\leq \frac{M+1}{M}\ell(\Delta)<2\ell(\Delta)$, and since $\ell(\Delta')/\ell(\Delta)$ is a power of $c_{0}<1/1000$, we must have $\ell(\Delta')\leq \ell(\Delta)$.
\item[$(\ddag)$] The above estimate also implies $\zeta_{\Delta}\in B(\zeta_{\Delta'},M\ell(\Delta))=B(\zeta_{\Delta'},5Mc_{0}^{n})$, and so the collection $\{c_{1}B_{\Delta}:\Delta\in \cM_{n}(\Delta')\}$ form a disjoint family of balls of radii $5c_{0}^{n}$ contained in $B(\zeta_{\Delta'},5(M+1)c_{0}^{n})$. Moreover, using the doubling property of $\mu$, and since $\ell(\Delta)=5c_{0}^{n}$,
\[\mu(c_{1}B_{\Delta})
\gec_{\mu,M,c_{1},c_{0}} \mu(4M B_{\Delta}) 
\geq \mu(B(\zeta_{\Delta'},5(M+1)c_{0}^{n})),\]
and thus we know that 
\begin{multline*}
\# \cM_{n}(\Delta') \mu(B(\zeta_{\Delta'},5(M+1)c_{0}^{n}))
\\ \lec_{M,c_{1},c_{0},\mu} \sum_{\Delta\in \cM_{n}(\Delta')}\mu(c_{1}B_{\Delta})
\leq \mu(B(\zeta_{\Delta'},5(M+1)c_{0}^{n}))
\end{multline*}
which implies
\begin{equation}
\label{e:numM}
\# \cM_{n}(\Delta')\lec_{\mu,M,c_{1},c_{0}} 1.
\end{equation} 
\end{enumerate}
Hence, for $\Delta_{1}\subseteq \Delta_{0}$, and $\beta>\beta_{0}$
\begin{align*}
\sum_{\Delta \subseteq \Delta_{1} \atop \Delta\cap E\neq\emptyset}\lambda(\Delta )\mu(\Delta )
& \stackrel{\eqn{b0}}{\lec}_{\mu,M,\beta} \sum_{\Delta \subseteq \Delta_{1} \atop \Delta\cap E\neq\emptyset}\sum_{\Delta' \subseteq MB_{\Delta} \atop \Delta'\in \cW(E^{c})} \frac{\mu(\Delta' )}{\mu(\Delta)} \ps{\frac{\ell(\Delta')}{\ell(\Delta)}}^{\beta-\beta_{0}}\mu(\Delta )\\
& = \sum_{\Delta' \in \cW(E^{c}) \atop \Delta'\subseteq MB_{\Delta_{1}}}\mu(\Delta' )\sum_{\Delta' \subseteq MB_{\Delta}, \Delta\cap E\neq\emptyset \atop \Delta \subseteq \Delta_{1}} \ps{\frac{\ell(\Delta')}{\ell(\Delta)}}^{\beta-\beta_{0}}\\
& \stackrel{(\dag)}{=} \sum_{\Delta' \in \cW(E^{c}) \atop \Delta'\subseteq MB_{\Delta_{1}}}\mu(\Delta' )\sum_{c_{0}^{n}\geq \ell(\Delta')} \sum_{\Delta\in \cM_{n}(\Delta')} \ps{\frac{\ell(\Delta')}{\ell(\Delta)}}^{\beta-\beta_{0}}\\
&  \stackrel{\eqn{numM}}{\lec}_{\mu,M,c_{1},c_{0}} \sum_{\Delta' \in \cW(E^{c}) \atop \Delta'\subseteq MB_{\Delta_{1}}}\mu(\Delta' ) \sum_{n\geq 0}c_{0}^{n(\beta-\beta_{0})}\\
& \lec_{\mu,\beta} \sum_{\Delta' \in W(E^{c}) \atop \Delta'\subseteq MB_{\Delta_{1}}}\mu(\Delta' )
\leq \mu(MB_{\Delta_{1}})
\lec_{\mu,M,c_{1},c_{0}} \mu(\Delta_{1}).
\end{align*}
\end{proof}

As a corollary, we have the following:

\begin{corollary}
With the assumptions of \Lemma{porous}, let $0<\delta<1<M<M_{0}/2$ and set
\[\cP_{M,\delta}=\{\Delta: \Delta\cap E\neq\emptyset, \exists \; \xi\in MB_{\Delta} \mbox{ such that }\dist(\xi,E)\geq \delta\ell(\Delta)\}.\]
Then there is $C_{1}=C_{1}(M,\delta,C_{\mu})>0$ so that, for all $\Delta'\subseteq \Delta_{0}$,
\begin{equation}
\sum_{\Delta\subseteq \Delta' \atop \Delta\in \cP_{M,\delta}} \mu(\Delta)\leq C_{1} \mu(\Delta').
\label{e:sumP}
\end{equation}
\label{c:porous}
\end{corollary}

\begin{proof}
If $\Delta\in \cP_{M,\delta}$, then there is $\xi\in MB_{\Delta}$ so that $\dist(\xi,E)\geq \delta\ell(\Delta)$, and so $B(\xi,\delta)\subseteq (M+\delta)B_{\Delta}\backslash E\subseteq 2MB_{\Delta}\backslash E$. Let $\Delta'$ be the maximal cube containing $\xi$ so that $2MB_{\Delta'}\cap E=\emptyset$. Then $\ell(\Delta')\leq \ell(\Delta)$ since, if $\zeta\in \Delta\cap E$,
\begin{multline*}
 2M\ell(\Delta')-\ell(\Delta)\leq \dist(\zeta_{\Delta'},E)-\ell(\Delta)
\leq |\zeta_{\Delta'}-\zeta|-|\zeta-\zeta_{\Delta}|\\
\leq |\zeta_{\Delta'}-\zeta_{\Delta}|\leq |\zeta_{\Delta'}-\xi|+|\xi-\zeta_{\Delta}| <\ell(\Delta')+M\ell(\Delta).
\end{multline*}
This and the fact that $\xi\in \Delta'\cap MB_{\Delta}$ imply $\Delta'\subseteq 2MB_{\Delta}$, and thus $\lambda_{2M,\beta}(\Delta)\geq (\ell(\Delta')/\ell(\Delta))^{\beta}$ where we set $\beta=2\log_{2}C_{\mu}$. 
Note that $\ell(\Delta')\geq \frac{\delta c_{0}}{2M}\ell(\Delta)$, since otherwise if $\Delta''$ is the parent of $\Delta'$, then $\xi\in \Delta''$ and so
\[MB_{\Delta''}\subseteq B(\xi,2M\ell(\Delta''))=B(\xi,2Mc_{0}\ell(\Delta'))\subseteq B(\xi,\delta\ell(\Delta))\subseteq E^{c},\]
but we know that, since $\Delta'$ is maximal, $MB_{\Delta''}\cap E\neq\emptyset$, and we get a contradiction. Thus, we have shown $\lambda_{2M,\beta(\Delta)}(\Delta)\geq \ps{\frac{\delta c_{0}}{2M}}^{\beta}$ whenever $\Delta\in \cP_{M,\delta}$, and the previous lemma implies \eqn{sumP}.
\end{proof}

\begin{lemma}
With the assumptions of \Lemma{porous}, and supposing $E\subseteq c_{0}B_{\Delta_{0}}\subseteq \Delta_{0}\in \cD$ is a Borel set satisfying $\mu(E)/\mu(\Delta_{0})\geq \rho>0$, we have that for all $\delta,\tau>0$, there are $t_{0},N>0$ (depending on $\delta,C_{\mu},M$, and $\rho$) such that for $t\in(0,t_{0})$, we can find a collection $T$ of cubes in $\Delta_{0}$, and a compact set $E'\subseteq E$ so that the following are true.
\begin{enumerate}
\item $\mu(E')\geq (1-\tau)\mu(E)$.
\item  If $\xi\in \Delta\cap E'$ for some $\Delta\in T$, then $\xi\in (1-t)\Delta$.
\item If $\Delta\subseteq \Delta_{0}$ and $\Delta\cap E'\neq \emptyset$, then either $\Delta\in T$ or, for every $\xi\in MB_{\Delta}$, $\dist(\xi,E)<\delta\ell(\Delta)$. 
\item For all $\Delta'\in \cD$,
\begin{equation}
\sum_{\Delta\subseteq \Delta' \atop \Delta\in T} \mu(\Delta)\lec_{\mu,\tau,\delta,M} \mu(\Delta').
\label{e:wpacking}
\end{equation}
\item Finally, we also have that, for every $\xi\in E'$, $\xi$ is contained in at most $N$ cubes from $T$.
\end{enumerate}
\label{l:E'}
\end{lemma}

\begin{proof}
Let $\cP=\cP_{M,\delta}$ be from \Corollary{porous}. For $\Delta\subseteq \Delta_{0}$, let $k(\Delta)$ denote the number of cubes in $\cP$ properly containing $\Delta$ (so $k(\Delta_{0})=0$). For $N>0$, by Chebyshev's inequality,
\[
\mu\ps{ \bigcup_{\Delta\subseteq \Delta_{0} \atop k(\Delta)\geq N} \Delta}
\leq \frac{1}{N}\sum_{\Delta\in \cP} \mu(\Delta) \stackrel{\eqn{sumP}}{\leq} \frac{C_{1}}{N}\mu(\Delta_{0}).\]
Thus, if $\tau\in (0,1)$, $N>\frac{2C_{1}}{\tau\rho}$, and
\[E_{N}:=E\backslash \bigcup_{\Delta\subseteq \Delta_{0} \atop k(\Delta)\geq N} \Delta\]
then
\[\mu(E_{N})\geq (1-\tau/2)\mu(E).\]
Set 
\[T=\{\Delta \subseteq \Delta_{0}:\Delta\cap E_{N}\neq\emptyset\}\cap \cP.\] 
By \Lemma{local-dub} and \eqn{sumP}, if $t<t_{0}:=\ps{\frac{\tau\rho}{2C_{1}t_{0}}}^{1/\alpha}$, then
\begin{equation}
\sum_{\Delta \in \cP} \mu(\Delta \backslash (1-t) \Delta)
\leq t_{0} t^{\alpha} \sum_{\Delta \in T} \mu(\Delta)
\leq C_{1} t_{0}t^{\alpha} \mu(\Delta_{0})<\frac{\tau}{2}\mu(E).
\label{e:w-ring}
\end{equation}
Thus, if
\begin{equation}
E'=E_{N}\cap\ps{ \bigcup_{\Delta\in T}\Delta\backslash (1-t) \Delta}^{c}
\label{e:E'}
\end{equation}
then
\[\mu(E')>(1-\tau)\mu(E).\]
Note that $E'\subseteq E_{N}$ guarantees $\xi$ is in at most $N$ many cubes from $T$. Finally, by replacing $E'$ with a compact subset if necessary so that it still satisfies the previous inequality, we may assume $E'$ is compact.
\end{proof}

\begin{remark}
It is this set of lemmas concerning porosity where the doubling property for harmonic measure (and hence the NTA assumption) plays the most critical role in our work. By work in \cite{BL04}, for example, one can generalize the results of \cite{DJ90} and \Theorem{DJ} to domains satisfying only an interior corkscrew condition, but whose boundary is Ahlfors regular and has ``uniform interior pieces of Lipschitz graphs," {\it a priori}. In this setting, harmonic measure isn't necessarily doubling, and so one only obtains a ``weak" $A_{\infty}$-condition or ``weak" reverse H\"older inequality (which implies the stronger $A_{\infty}$-condition if $\omega$ happens to be doubling). Thus, one could perhaps generalize our results in this way via constructing $\Omega_{E'}^{\pm}$ satisfying interior corkscrew conditions and using the comparison principle; however, we also use the doubling property to construct these ideal subsets $E'$ that guarantee that our domains $\Omega_{E'}^{\pm}$ have Ahlfors regular boundaries. It's because of this that a generalization is even less immediate.
\label{r:generalize}
\end{remark}

\section{The sets $\Omega_{E}^{\pm}$}
\label{s:omegaE}

We will use a pretty general method for constructing sub and super NTA domains that intersect a prescribed portion of the boundary, and later prove that, given a clever choice of subset $E'\subseteq E$ (where $E$ is as in \Theorem{main}), the sub and super NTA domains containing $E'$ in their boundaries have the desired properties. The constructions of the subdomains are common knowledge (see \cite{HM12} or \cite{JerisonKenig}, for example), but the existence of superdomains is not, to the author's knowledge, and so we include a construction in the appendix.

\begin{lemma}
Let $\Omega\subseteq\bR^{d+1}$ be a $C$-NTA (or $C$-uniform) domain and let $E\subseteq B(\xi_{0},r_{0})\cap \d\Omega$ be compact where $\xi_{0}\in \d\Omega$ and $r_{0}\in(0,\diam\d\Omega)$. Set $C_{0}>0$ and 
\[\cC_{E}^{-}=\{Q\in \cW(\Omega): C_{0}Q\cap E\neq\emptyset, \ell(Q)\leq  r_{0}\}.\]
For $Q_{1},Q_{2}\in \cW(\Omega)$, let $P_{Q_{1},Q_{2}}$ be a shortest path of adjacent dyadic Whitney cubes connecting $Q_{1}$ to $Q_{2}$ (which also includes $Q_{1}$ and $Q_{2}$). For some constant $\tilde{C}>0$, set
\[\widetilde{\cC_{E}}^{-}=\{Q:Q\in P_{Q_{1},Q_{2}}\mbox{ for some }Q_{1},Q_{2}\in \cC_{E}^{-} \mbox{ with } d_{\Omega}(Q_{1},Q_{2})\leq \tilde{C}\}.\]
For $\lambda>1$, set
\[\Omega_{E}^{-}=\ps{\bigcup_{Q\in \widetilde{\cC_{E}}^{-}}\lambda Q}^{\circ}.\]
Then for $C_{0}$ and $\tilde{C}$ large enough and $\lambda>1$ close enough to $1$(each depending only on $C$ and $d$), $\Omega_{E}^{-}$ is a $C^{-}$-NTA (or $C^{-}$-uniform) domain contained in $B(\xi_{0},C^{-}r_{0})$ for some $C^{-}=C^{-}(d,C_{0},\lambda,C)$ and $\diam \d\Omega_{E}^{-}\geq r_{0}/C^{-}$ . Moreover, $\d\Omega_{E}^{-}\cap \d\Omega= E$. 
\label{l:omegaE-}
\end{lemma}

As mentioned earlier, we will omit the proof of \Lemma{omegaE-}, but for an idea of the construction, see \cite[Lemma 6.3]{JerisonKenig} or \cite[Lemma  3.61]{HM12} for example. 

\begin{lemma}
Let $\Omega$ be a $C$-NTA (or $C$-uniform) domain, and $E\subseteq \d \Omega\cap B(\xi_{0},r_{0})$ be compact where $\xi_{0}\in \d\Omega$ and $r_{0}\in (0,\diam \d\Omega)$. Let $K\geq 3,\lambda>1$ and set
\[ \cC_{E}^{+}=\{Q\in \cW_{K}(\bR^{d+1}\backslash E): Q\cap \d\Omega\neq\emptyset\}.\]
Define
\[\Omega_{E}^{+}=\Omega\cup \bigcup_{Q\in \cC_{E}^{+}} (\lambda Q)^{\circ}.\]
Then, for $\lambda>1$ close enough to $1$ (depending on $C$ and $d$) and $K$ large enough (depending on $d,\lambda$ and $C$), there is $C^{+}=C^{+}(d,C,K)$ so that $\Omega_{E}^{+}$ is a $C^{+}$-NTA (or $C^{+}$-uniform) domain. Moreover, $\d\Omega_{E}^{+}\cap \d\Omega={E}$. If $\Omega$ is $C$-NTA, then $\diam \d\Omega_{E}^{+}\sim_{C^{+}} \diam \d\Omega$. 
\label{l:omegaE+}
\end{lemma}

For a proof, see the appendix.

\begin{lemma}
Let  $C'>0$, $\Omega$, $E$,  and $\Omega_{E}^{\pm}$ be as in \Lemma{omegaE-} or \Lemma{omegaE+}. Suppose also that there is an $L$-bi-Lipschitz injection $g:E\rightarrow Z$ where $Z$ is a metric space satisfying \eqn{almost-ahlfors} and that for all $\xi\in E$ and $r>0$, we have 
\begin{equation}
\sum_{Q\in \d\cC_{E}^{\pm}(\xi,r)}\ell(Q)^{d}\leq C'r^{d}
\label{e:dCEsum}
\end{equation}
where
\begin{multline}
\d{\cC_{E}^{-}}(\xi,r)
=\{Q\in \cC_{E}^{-}: Q\cap B(\xi,r)\neq\emptyset,Q\sim Q'\mbox{ for some } \\
Q'\in \cW(\Omega)\backslash \cC_{E}^{-}\}
\end{multline}
and 
\begin{multline}
\d{\cC}_{E}^{+}(\xi,r)
=\{Q\in {\cC_{E}}^{+}: Q\cap B(\xi,r)\neq\emptyset, Q\sim Q'\mbox{ for some } \\
Q'\in \cW_{K}(\bR^{d+1}\backslash E)\backslash\cC_{E}^{+}\}.
\end{multline}
Then $\d\Omega_{E}^{\pm}$ is upper $A^{\pm}$-Ahlfors regular (with $A^{\pm}=A^{\pm}(C,C',A,L,d)$), meaning 
\begin{equation}
\cH^{d}(B(\xi,r)\cap \d\Omega_{E}^{\pm})\leq A^{\pm} r^{d} \mbox{ for all }\xi\in \d\Omega_{E}^{\pm}\mbox{ and }r>0.
\label{e:uregular}
\end{equation}
 If $\Omega$ is also $C$-NTA, then $\d\Omega_{E}^{\pm}$ is $A^{\pm}$-Ahlfors regular. 
 
\label{l:lqsum}

\end{lemma}

\begin{proof}

\Claim: For all $r>0$, if
\begin{multline}
\d\widehat{\cC}_{E}^{-}(\xi,r)
=\{Q\in \widetilde{\cC_{E}}^{-}:\lambda Q\cap B(\xi,r)\neq\emptyset, Q\sim Q'\mbox{ for some } \\
Q'\in \cW(\Omega)\backslash\widetilde{\cC_{E}}^{-}\}
\label{e:dhatce-}
\end{multline}
and
\begin{multline}
\d\widehat{\cC}_{E}^{+}(\xi,r)
=\{Q\in \d{\cC_{E}}^{+}:\lambda Q\cap B(\xi,r)\neq\emptyset, Q\sim Q'\mbox{ for some } \\
Q'\in \cW_{K}(\bR^{d+1}\backslash E)\backslash\d\cC_{E}^{+}\}
\label{e:dhatce+}
\end{multline}
then
\begin{equation}
\sum_{Q\in \d\widetilde{\cC}_{E}^{\pm}(\xi,r)}\ell(Q)^{d}\lec_{C',d,K,\tilde{C},C} r^{d}.
\label{e:hatsum}
\end{equation}

We first focus on $\d\cC_{E}^{+}(\xi,r)$. If $Q\in \d\widehat{\cC}_{E}^{+}(\xi,r)$, then $\lambda Q\cap B(\xi,r)\neq\emptyset$, and by \Definition{Whitney}, 
\begin{equation}
\ell(Q) \leq  r/\ps{\frac{K-1}{2}-\sqrt{d+1}\frac{\lambda-1}{2}}<4r/(K-1)
\label{e:4r/k-1}
\end{equation}
for $\lambda>1$ close enough to $1$, and so $\diam \lambda Q\leq 4\lambda\sqrt{d+1}r/(K-1)$. Hence $Q\subseteq B(\xi,(4\lambda \sqrt{d+1}/(K-1)+1)r)$, and so
\[\sum_{Q\in \d\widehat{\cC}_{E}^{+}(\xi,r)}\ell(Q)^{d}
\leq \sum_{Q\in \d \cC_{E}^{+}(\xi,(4\lambda \sqrt{d+1}/(K-1)+1)r)}\ell(Q)^{d} \stackrel{\eqn{dCEsum}}{\lec}_{K,d} r^{d}\]
which proves the claim in this case.

In the case of $\d\cC_{E}^{-}(\xi,r)$, if $Q\in \d \widehat{\cC}_{E}^{-}(\xi,r)$, \eqn{4r/k-1} still holds with $K=3$. Moreover, there is a chain of Whitney cubes of length $\tilde{C}$ of Whitney cubes from $Q$ to a cube $Q'\in \d \cC_{E}^{-}$, each cube in the chain having diameter comparable to $\ell(Q)$ (with constants depending on $d$ and $\tilde{C}$), so in particular, $\dist(Q,Q')\lec_{\tilde{C},d} \ell(Q)\leq 2r$ and $\ell(Q')\sim_{\tilde{C},d}\ell(Q)$, and so
\[\dist(\xi,Q')\leq \dist(\xi,Q)+\diam Q+\dist(Q,Q')\lec_{\tilde{C},d} r.\]
Thus, there is $C''$ depending on $d$ and $\tilde{C}$ so that $Q'\subseteq B(\xi,C''r)$. Also, to each $R\in \d\cC_{E}^{-}$, there are at most $N=N(\tilde{C},d)$ cubes $Q\in \d\widehat{\cC}_{E}^{-}$ with $Q'=R$. Thus,
\[\sum_{Q\in \d\widehat{\cC}_{E}^{-}(\xi,r)}\ell(Q)^{d}
\sim_{\tilde{C},d}\sum_{Q\in \d\widehat{\cC}_{E}^{-}(\xi,r)}\ell(Q')^{d}
\lec_{\tilde{C},d} \sum_{R\in \d \cC_{E}^{-}(\xi,C''r)}\ell(R)^{d} \stackrel{\eqn{dCEsum}}{\lec}_{\tilde{C},d} r^{d}.\]
Thus we've finished the claim.

Now we will prove \eqn{uregular}.

\begin{enumerate}
\item Suppose $\dist(\xi,E)\geq 2r$. Let $Q\in \d \widehat{\cC}_{E}^{\pm}(\xi,r)$ and $y_{Q}\in B(\xi,r)\cap \d \lambda Q$. Then $Q\in \cW(\Omega)$ or $Q\in \cC_{E}^{+}\subseteq \cW_{K}(\bR^{d+1}\backslash E)$ (depending on whether we're considering $\Omega_{E}^{-}$ or $\Omega_{E}^{+}$; if the former, we set $K=3$), so by \Definition{Whitney}, for $\lambda>1$ small enough (recall $K\geq 3$)
\[  \ell(Q)\leq \frac{\dist(Q,E)}{\frac{K-1}{2}-\sqrt{d+1}\frac{\lambda -1}{2}}\leq 2\dist(Q,E) \leq 2r.\]
Thus, since $\lambda Q \cap B(\xi,r)\neq\emptyset$,
\[\lambda Q\subseteq B(\xi,r+\diam \lambda Q)\subseteq B(\xi,(1+2\lambda \sqrt{d+1})r).\]
Moreover,
\[\ell(Q)\geq \frac{\dist(y_{Q},E)}{(1+K+(\lambda-1)/2)} \gec_{K} \dist(\xi,E)-|y_{Q}-\xi|> r.\]
Thus, $B(\xi,r)\cap \d\Omega_{E}^{\pm}$ is in the union of the boundaries of finitely many cubes of the form $\lambda Q$ where the $Q$ have diameters comparable to $r$ and is contained in a ball of radius comparable to $r$; this implies $\cH^{d}(B(\xi,r)\cap \d \Omega_{E}^{\pm})\lec_{d,\lambda} r^{d}$ (where the implied constant depends on $K$ in the case of $\Omega_{E}^{+}$).
\item  If $\dist(\xi,E)< 2r$, let $\xi'\in E$ be such that $|\xi'-\xi|< 2r$. Then

\begin{align*}
\cH^{d}(\d \Omega_{E}^{\pm}\cap B(\xi,r))
& \leq \cH^{d}(\d \Omega_{E}^{\pm}\cap B(\xi',2r))\\
& \leq \sum_{Q\in \d\widehat{\cC}_{E}^{\pm}(\xi',2r)}\cH^{d}(\d \lambda Q)+\cH^{d}(E\cap B(\xi',3r))\\
& \lec_{\lambda,d} \sum_{Q\in \d\widehat{\cC}_{E}^{\pm}(\xi',2r)}\cH^{d}(\d \lambda Q)+ L^{d}\cH_{Z}^{d}(g(E\cap B(\xi',2r)))\\
& \stackrel{\eqn{hatsum}}{\lec}_{\lambda,d} r^{d}+ L^{d}\cH_{Z}^{d}(B_{Z}(\xi',2Lr)) 
\stackrel{\eqn{almost-ahlfors}}{\lec}_{A,L,d}  r^{d}\end{align*}
and this proves \eqn{uregular}. Note that \eqn{almost-ahlfors} is given only for radii at most $r_{0}$, but since $E\subseteq B(\xi_{0},r_{0})$, it also holds for all $r>0$ with perhaps a slightly larger constant.
\end{enumerate}

This proves the lemma for the case of $C$-uniform $\Omega$. If $\Omega$ is $C$-NTA, so are $\Omega_{E}^{\pm}$ and it is well known that $C$-NTA domains are lower regular (that is, the lower bound in \eqn{regular}) with constant depending on $C$ and $d$. To see this, let $B(x,r/C)\subseteq B(\xi,r)\cap \Omega_{E}^{\pm}$ and $B(y,r/C)\subseteq B(\xi,r)\backslash \Omega_{E}^{\pm}$, let $P_{x}$ and $P_{y}$ be two parallel $d$-planes passing through $x$ and $y$ respectively, and let $D_{x}=B(x,r/C)\cap P_{x}$ and $D_{y}\cap P_{y}\cap B(y,r/C)$. Note that each segment perpendicular to $P_{x}$ and passing from $D_{x}$ to $D_{y}$ must intersect $\d\Omega_{E}^{\pm}$, and thus if $\pi$ is the orthogonal projection onto $P_{x}$,
\[\cH^{d}(\d\Omega_{E}^{\pm}\cap B(\xi,r))\geq \cH^{d}(\pi(\d\Omega_{E}^{\pm} \cap B(\xi,r)))\geq \cH^{d}(D_{x})\gec_{C,d} r^{d}.\]

\end{proof}

Thus, the main challenge in proving \Theorem{uniform-main} is to show how our assumptions imply \eqn{dCEsum} holds for $E$, or in the case of \Theorem{main}, to show \eqn{dCEsum} holds for some special subset $E'$.

\section{The proof of \Theorem{main}}
\label{s:main}

We now apply the results of the previous sections to prove \Theorem{main}. We state here our standing assumptions that will hold throughout this section:\\

{\bf Standing assumptions for this section:} We will assume $\Omega$ is a $C$-NTA domain, $E\subseteq \d\Omega\cap B(\xi_{0},r_{0})$, $g:E\rightarrow Z$ is a $L$-bi-Lipschitz injection into a metric space $Z$ satisfying \eqn{almost-ahlfors}, $r_{0}\in (0,\diam \d\Omega)$, $\xi_{0}\in \d\Omega$, $B(z_{0},r_{0}/C)\subseteq \Omega$, $\omega=\omega_{\Omega}^{z_{0}}$, and $\omega(E)/\omega(B(\xi_{0},r_{0}))\geq \rho>0$. We will also assume that $\cD$ are the ``cubes" for $\Sigma=\d\Omega$ with $c_{0}<c_{1}/4$ fixed, and without loss of generality, that $B(\xi_{0},r_{0})=B(\zeta_{\Delta_{0}},c_{1}\ell(\Delta_{0}))$ (we can do this by rescaling $\Omega$ and by choosing the maximal nets in \Theorem{Christ} to include $\xi_{0}$). We will also let $B_{\Delta}=B(\zeta_{\Delta},\ell(\Delta))$ denote a {\it Euclidean ball}, not a ball with respect to the relative topology of $\d\Omega$, though we still have $\d\Omega\cap c_{1} B_{\Delta}\subseteq \Delta\subseteq B_{\Delta}\cap \d\Omega$ for $\Delta\in\cD$, so in particular, $E\subseteq c_{0}B_{\Delta_{0}}\cap \Sigma\subseteq \Delta_{0}$. Let $M>0$ be large and $\delta>0$ to be determined later.

Note that by \eqn{wdub}, we can apply \Lemma{porous} and  \Lemma{E'} with $\mu=\omega$, $\tau>0$, some numbers $\delta,M>0$ to be chosen later, and $\rho \frac{\omega(B(\xi_{0},r_{0}))}{\omega(\Delta_{0})}$ in place of $\rho$. From this lemma, we obtain the quantities $t_{0}$ and $N$, and for $t\in(0,t_{0})$, we get a compact set $E'$ and a collection of cubes $T$, where for now we pick $t\in(0,t_{0})$ small enough (depending on $c_{1}$) so that
\begin{equation}
(c_{1}/2)B_{\Delta}\cap \d\Omega \subseteq \Delta\backslash (1-t)\Delta \mbox{ for }\Delta\in \cD
\label{e:bigcenter}
\end{equation}

In \Theorem{main}, the last conclusion follows from the penultimate one, and that one follows from the first five conclusions and \Lemma{afineq}, so we need only prove those. Let $\Omega_{E'}^{\pm}$ be the NTA domains from Lemmas \ref{l:omegaE-} and \ref{l:omegaE+}, where we will pick $C_{0}$ in the course of the proof sufficiently large. These and the set $E'$ already satisfy conclusions (1) through (4), so we only have to demonstrate that they are Ahlfors regular. By \Lemma{lqsum}, it suffices to show \eqn{dCEsum}.

With all these reductions and assumptions in place, \Theorem{main} will now follow from the following lemma

\begin{lemma}
Fix $\tau\in (0,1)$, let $E'\subseteq E$ be the set from \Lemma{E'} for our choice of $\tau$, $\rho \frac{\omega(B(\xi_{0},r_{0}))}{\omega(\Delta_{0})}$ in place of $\rho$, and some $M$ and $\delta$ and $\mu=\omega$, and $\Omega_{E'}^{\pm}$ be as in \Lemma{omegaE-} or \Lemma{omegaE+} applied to the set $E'$.  Then
\begin{equation}
\sum_{Q\in \d\cC_{E'}^{\pm}(\xi,r)}\ell(Q)^{d}\lec r^{d} \mbox{ for all }\xi\in E'\mbox{ and }r>0.
\label{e:QindC}
\end{equation}
\label{l:lqsum2}
\end{lemma}

\begin{proof}[Proof of \Lemma{lqsum2}]

The sum in \eqn{QindC} will be controlled using two different bookkeeping lemmas.

\begin{lemma}
Let $E'\subseteq E\subseteq B(\xi_{0},r_{0})\subseteq \bR^{d+1}$, and $r\in (0,3r_{0})$, and $\cC$ be any collection of disjoint cubes $Q\subseteq\bR^{d+1}$. Suppose $y_{Q}\in E\backslash \cnj{E'}$ are points in $E$ such that 
\begin{enumerate}
\item $\dist(y_{Q},Q)\lec \dist(y_{Q},E')\sim \dist(Q,E')\sim \ell(Q)$,
\item for all $Q\in \cC$, $\dist(Q,E') < r$, 
\item there is an $L$-bi-Lipschitz injection $g:E\rightarrow Z$ into a metric space $Z$ satisfying \eqn{almost-ahlfors}.
\end{enumerate}
Then $\sum_{Q\in \cC}\ell(Q)^{d}\lec r^{d}$ (with constant depending on $d$, $L$, the constants in \eqn{almost-ahlfors}, and all implied constants).
\label{l:Zlemma}
\end{lemma}

\begin{proof}
Let $\cD_{Z}$ denote the ``cubes" for $Z$ and $F=g(E')$. Let $\cW(Z\backslash F)\subseteq\cD_{Z}$ denote the collection of maximal cubes $\Delta$ for which $3B_{\Delta}\cap F=\emptyset$. One can show $\ell(\Delta)\sim \dist(\xi,F)$ for all $\xi \in \Delta\in \cW(Z\backslash F)$. Let $\Delta_{Q}\in \cW(Z\backslash F)$ contain $g(y_{Q})$; we know such a cube exists since \[\dist(g(y_{Q}),F)\geq L^{-1}\dist(y_{Q},E')>0.\]
Note that 
\[\ell(Q)\sim \dist(y_{Q},E')\sim \dist(g(y_{Q}),F)\sim \ell(\Delta_{Q}).\]
\Claim There is $N_{0}=N_{0}(d,L)$ so that at most $N_{0}$ many $Q\in \cC$ can satisfy $\Delta_{Q}=\Delta$ for some given $\Delta\in \cW(Z\backslash F)$. To see this, note
\[\dist(y_{Q},Q)\lec \dist(y_{Q},E')\sim \ell(\Delta_{Q})=\ell(\Delta),\]
and if $\Delta_{Q}=\Delta_{Q'}=\Delta$, then $|y_{Q}-y_{Q}'|\leq L|g(y_{Q})-g(y_{Q'})|\leq \diam\Delta\sim \ell(\Delta)$. Thus, all cubes $Q$ for which $\Delta_{Q}=\Delta$ are contained in a ball of radius comparable to $\ell(\Delta)$ and have side lengths comparable to $\ell(\Delta)$. This proves the claim.

Since each $Q\in \cC$ intersects $B(\xi,r)$,
\[\ell(\Delta_{Q})\sim \dist(y_{Q},E') \sim \ell(Q) \sim \dist(Q,E')\leq r\]
and if we fix a $Q_{0}\in \cC$,
\begin{multline*}
\dist(g(y_{Q_{0}}),\Delta_{Q}) 
\leq |g(y_{Q_{0}})-g(y_{Q})|\lec_{L} |y_{Q_{0}}-y_{Q}|\\
 \leq \dist(y_{Q_{0}},Q_{0})+\diam Q_{0}+\dist(Q_{0},Q)+\diam Q+\dist(y_{Q},Q)\\
 \lec \ell(Q_{0})+\ell(Q_{0})+2r+\ell(Q)+\ell(Q)\lec r
\end{multline*}
thus, all the $\Delta_{Q}$ are contained in a ball $B=B_{Z}(g(\xi), C'r)$ for some $C'$ depending on $L,d$ and the implied constants. Moreover, since $\ell(\Delta_{Q})\sim \ell(Q)\leq r<3r_{0}$, there is $\theta<c_{1}$ so that $\theta \ell(\Delta_{Q})<r_{0}$, ($\theta<c_{1}$ guarantees $\theta B_{\Delta}\subseteq \Delta$). Thus, if $C'r<r_{0}$, we can apply part (3) to get
\begin{align*}
\sum_{Q\in\cC}\ell(Q)^{d}
& \lec_{N_{0}} \sum_{\Delta\in \cW(Z\backslash F) \atop \Delta \subseteq B}\ell(\Delta)^{d}
\lec_{A}  \sum_{\Delta\in \cW(Z\backslash F) \atop \Delta \subseteq B} \cH_{Z}^{d}(\theta B_{\Delta}) \\
& \leq  \sum_{\Delta\in \cW(Z\backslash F) \atop \Delta \subseteq B} \cH_{Z}^{d}(\Delta) 
\leq \cH_{Z}^{d}(B) \lec_{A} r^{d}.
\end{align*}
If $C'r\geq r_{0}$, we can cover $E$ with a bounded number (depending on $d$ and $L$) of balls $B_{i}$ centered on $E$ with radii less than $\frac{r_{0}}{2L}$, and thus we can cover $g(E)$ with a finite number of balls $B_{i}'$ centered on it of radii less than $r_{0}/2$. For $\theta$ small enough, $\theta B_{\Delta_{Q}}\subseteq \bigcup 2B_{i}'$, and thus
\begin{align*}
\sum_{Q\in\cC}\ell(Q)^{d}
& \lec_{N_{0}} \sum_{\Delta\in \cW(Z\backslash F) \atop \theta B_{\Delta} \subseteq \bigcup B_{i}'}\ell(\Delta)^{d}
\lec_{A}  \sum_{\Delta\in \cW(Z\backslash F) \atop \theta B_{\Delta} \subseteq \bigcup B_{i}'} \cH_{Z}^{d}(\theta B_{\Delta}) \\
& \leq \cH^{d}_{Z}\ps{\bigcup B_{i}'}\lec r_{0}^{d}\lec_{A} r^{d}.
\end{align*}
\end{proof}

\begin{lemma}
Let $T$ be the collection of cubes from \Lemma{E'} and $T_{\Delta'}$ be those cubes in $T$ contained in $\Delta'$ that intersect $E'$. Then for all $\Delta'\subseteq \Delta_{0}$,
\begin{equation}
\sum_{\Delta\in T_{\Delta'}} \ell(\Delta)^{d}\lec_{L,d,C,\tau} \ell(\Delta')^{d}.
\label{e:better}
\end{equation}
\label{l:better}
\end{lemma}

\begin{proof}
First note that, if $\Delta\subseteq \Delta_{0}$ and $\Delta\cap E'\neq \emptyset$, then there is $\xi_{\Delta}\in (1-t)\Delta \cap E'$ by \Lemma{E'}. Hence, the collection $\{B(\xi_{\Delta},t\ell(\Delta)):\Delta\subseteq \Delta_{0},\Delta\in \cD_{n}, \Delta\cap E'\neq\emptyset\}$ is a disjoint family of balls with centers in $E'\subseteq Z$.

Let $B^{\Delta}:=B_{Z}(g(\xi_{\Delta}),\frac{t}{2L}\ell(\Delta) )\subseteq Z$. If $y\in B^{\Delta}\cap B^{\Delta'}$ for some $\Delta,\Delta'\in \cD_{n}\cap T$, then $|\xi_{\Delta}-\xi_{\Delta'}|< tc_{0}^{n}$, so that $B(\xi_{\Delta},t\ell(\Delta))$ and $B(\xi_{\Delta'},t\ell(\Delta'))$ intersect, giving a contradiction. Thus, for any $y\in Z$, there is at most one $\Delta\in \cD_{n}\cap T$ containing $y$, and thus there are at most $N$ cubes $\Delta$ from $T$ so that $y\in B^{\Delta}$. Also, note that $\Delta,\tilde{\Delta}\in T_{\Delta'}$ implies 
\[|g(\xi_{\Delta})-g(\xi_{\tilde{\Delta}})|\leq L|\xi_{\Delta}-\xi_{\tilde{\Delta}}|<2L\ell(\Delta')\]
and so all the $B_{\Delta}$ lie in $B_{Z}(g(\xi_{\tilde{\Delta}}),(2L+\frac{t}{2L})\ell(\Delta'))\subseteq \tilde{B}:=B_{Z}(g(\xi_{\tilde{\Delta}}),(2L+1)\ell(\Delta))$ for some fixed $\tilde{\Delta}\in T_{\Delta'}$. 

Recalling that $r_{0}=c_{1}\ell(\Delta_{0})$, pick $\theta=c_{1}/2$ so that for all $\Delta\subseteq\Delta_{0}$, $\theta\ell(\Delta)\leq \theta\ell(\Delta_{0})<r_{0}$. Then by \eqn{almost-ahlfors}, if $\zeta\in \Delta'\cap E'$ and $(2L+1)\ell(\Delta')<r_{0}$,
\begin{align*}
\sum_{\Delta\in T_{\Delta'}} & \ell(\Delta)^{d}
 \lec_{A,d,t,c_{1},L } \sum_{\Delta\in T_{\Delta'}} \cH_{Z}^{d}(\theta B^{\Delta})
 \leq \int_{Z}\sum_{\Delta\in T_{\Delta'}} \one_{B^{\Delta}}(x)d\cH_{Z}^{d}(x)\\
& \leq N\cH_{Z}^{d}\ps{ \bigcup_{\Delta\in T_{\Delta'}} B^{\Delta}}
 \leq N\cH_{Z}^{d}( \tilde{B})
 \lec_{A,d,N} \ell(\Delta')^{d}.
\end{align*}
Otherwise, if $(2L+1)\ell(\Delta)\geq r_{0}=c_{0}\ell(\Delta_{0})$, cover $\Delta_{0}$ with $N_{1}=N_{1}(d,L,c_{1})$ many cubes $\Delta_{j}$ with $(2L+1)\ell(\Delta_{j})<r_{0}$. Then
\[
\sum_{\Delta\in T_{\Delta'}}  \ell(\Delta)^{d}
\leq \sum_{j=1}^{N_{1}}\sum_{\Delta\in T_{\Delta_{j}}}  \ell(\Delta)^{d}
 \lec \sum_{j=1}^{N_{1}}\ell(\Delta_{j})^{d}\lec_{L,d,c_{1}} \ell(\Delta_{0})^{d}\lec_{L} \ell(\Delta').
\]
\end{proof}

\begin{lemma}
The inequality \eqn{QindC} holds for $\d\cC_{E'}^{-}(\xi,r)$.
\label{l:-case}
\end{lemma}

\begin{proof}
\Claim It suffices to show \eqn{QindC} in the case when $r\leq 3r_{0}$. To see this, observe that, if $r>3r_{0}$, then since $E'\subseteq E\subseteq B(\xi_{0},r_{0})$  any cube $Q\in \cC_{E'}^{-}\backslash \d\cC_{E'}^{-}(\xi,3r_{0})$ is at least $r_{0}$ away from $E'$. By construction, however, the $Q\in \cC_{E'}^{-}$ are chosen so that $\ell(Q)\leq r_{0}$. Since $C_{0}Q\cap E'\neq\emptyset$, we have that 
\[r_{0}\leq \dist(Q,E')\leq \dist(\xi,Q)\leq \diam C_{0}Q= C_{0}\sqrt{d+1} \ell(Q) \leq C_{0}\sqrt{d+1}r_{0}.\]
Thus all such $Q$ have sizes comparable to $r_{0}$. Moreover, all cubes in $ \cC_{E'}^{-}\backslash \cC_{E'}^{-}(\xi,r)$ lie in $B(\xi_{0},C^{-}r_{0})$ (since $\Omega_{E'}^{-}\subseteq B(\xi_{0},C^{-}r_{0})$) and thus there are boundedly many of them (depending only on $C_{0}$ and $d$). Hence, 
\[
\sum_{Q\in \d \cC_{E'}^{-}(\xi,r)}\ell(Q)^{d}
\leq \sum_{Q\in  \cC_{E'}^{-}\backslash \d\cC_{E}^{-}(\xi,3r_{0})}\ell(Q)^{d}+\sum_{Q\in \d \cC_{E'}^{-}(\xi,3r_{0})}\ell(Q)^{d}
\lec_{C_{0},d} r_{0}^{d}\lec_{d} r^{d}\]
which proves \eqn{QindC} if we assume \eqn{QindC} holds for $r\leq 3r_{0}$, and this proves the claim.

Now assume $r\leq 3r_{0}$. Set
\[\cD(\xi,r)= \{\Delta\in \cD(\Delta_{0}): \ell(\Delta)> r\geq c_{0}\ell(\Delta),\Delta\cap B(\xi,r)\cap E'\neq\emptyset\} .\]
Since $r\leq 3r_{0}=3c_{1}\ell(\Delta_{0})<\ell(\Delta_{0})$ and $E\subseteq \Delta_{0}$, this set is nonempty and covers $B(\xi,r)\cap E'$. Note that for $Q\in \d \cC_{E}^{-}(\xi,r)$, 
\[\ell(Q)\leq \dist(Q,\d\Omega)\leq \dist(\xi,Q)<r\leq 3r_{0} <\ell(\Delta_{0});\] 
this and the fact that $C_{0}Q\cap E'\neq\emptyset$ imply there is a maximal $\Delta(Q)\subseteq \Delta_{0}$ that intersects $C_{0}Q\cap E'$ and is such that $\ell(Q)\geq c_{0}\ell(\Delta(Q))$, so necessarily, $\Delta(Q)$ is contained in some cube in $\cD(\xi,r)$ (observe also that $\ell(\Delta(Q))\sim_{c_{0}} \ell(Q)$). With this in mind, and the fact that $\# \cD(\xi,r)\lec_{C} 1$, it will now suffice to show instead that
\begin{equation} \sum_{Q\in C_{E'}(\tilde{\Delta})} \ell(Q)^{d}\lec r^{d}\mbox{ for }\tilde{\Delta}\in \cD(\xi,r)
\label{e:CEDtwid}
\end{equation}
where
\[C_{E'}(\tilde{\Delta})=\{Q\in \d\cC_{E'}^{-}(\xi,r): \Delta(Q)\subseteq \tilde{\Delta}\}.\]

Split $C_{E'}(\tilde{\Delta})$ into sets 
\[T_{1}= \{Q\in C_{E'}(\tilde{\Delta}):\Delta_{Q}\in T\} , \; \;T_{2}=C_{E'}(\tilde{\Delta})\backslash T_{1}.\]

We first handle $T_{1}$. Observe that at most a bounded number of cubes $Q$ can have $\Delta(Q)=\Delta$ for a given $\Delta$ since $\dist(Q,\Delta(Q))\leq C_{0}\diam Q$ and $\ell(\Delta(Q))\sim \ell(Q)$. Also, since $\Delta(Q)\cap E'\neq\emptyset$ for all $Q\in C_{E'}(\tilde{\Delta})$, we know $T_{1}\subseteq T_{\tilde{\Delta}}$, (recall \Lemma{better} for this notation). Thus, we have
\begin{equation}
\sum_{Q\in T_{1}}\ell(Q)^{d}
\lec \sum_{\Delta\in T_{\tilde{\Delta}}}\ell(\Delta)^{d}\stackrel{\eqn{better}}{\lec} \ell(\tilde{\Delta})^{d}\lec r^{d}.
\label{e:T1}
\end{equation}

Next, assume $Q\in T_{2}$. Since $Q\in \d \cC_{E'}^{-}$, there is $Q'$ adjacent to $Q$ so that $C_{0}Q'\cap E' =\emptyset$. We now pick $C_{0}>0$ large enough (depending on $c_{0}$) so that there is $\Delta'\subseteq C_{0}Q'$ with $c_{0}\ell(\Delta')\leq \ell(Q')\leq \ell(\Delta')$ (see Figure \ref{f:qq}). Then 
\[\ell(\Delta')\sim_{c_{0}} \ell(Q')\sim_{d} \ell(Q)\sim_{c_{0}} \ell(\Delta(Q))\]
and so for $M>0$ large enough, $\Delta'\subseteq M B_{\Delta(Q)}$. Since $\Delta(Q)\not\in T$, if $\delta$ is small enough, then by \eqn{bigcenter} and \Lemma{E'}, there is $y_{Q}\in (1-t)\Delta'\cap E$. 

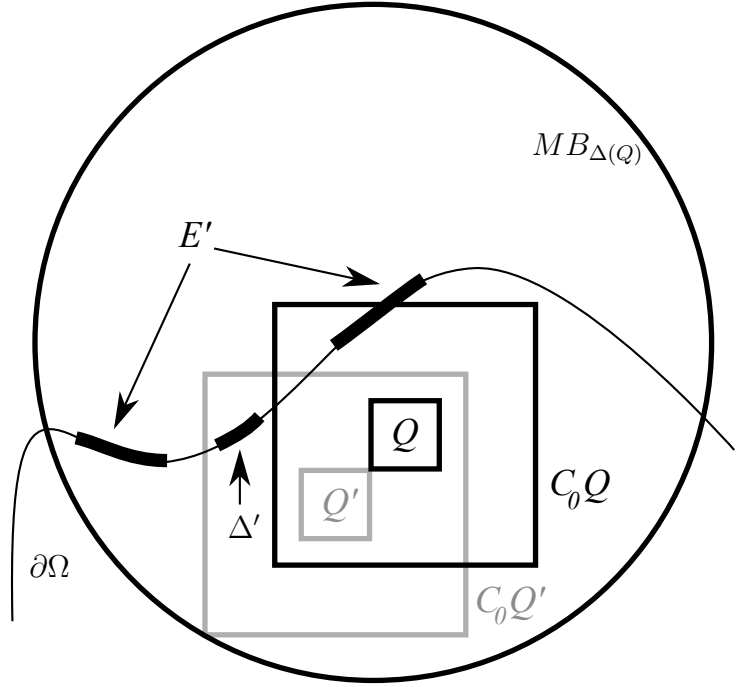
\begin{figure}[h]
\begin{picture}(0,0)(0,0)
\put(200,200){$MB_{\Delta(Q)}$}
\put(10,40){$\d \Omega$}
\put(85,55){$\Delta'$}
\end{picture}
\scalebox{1}{\includegraphics{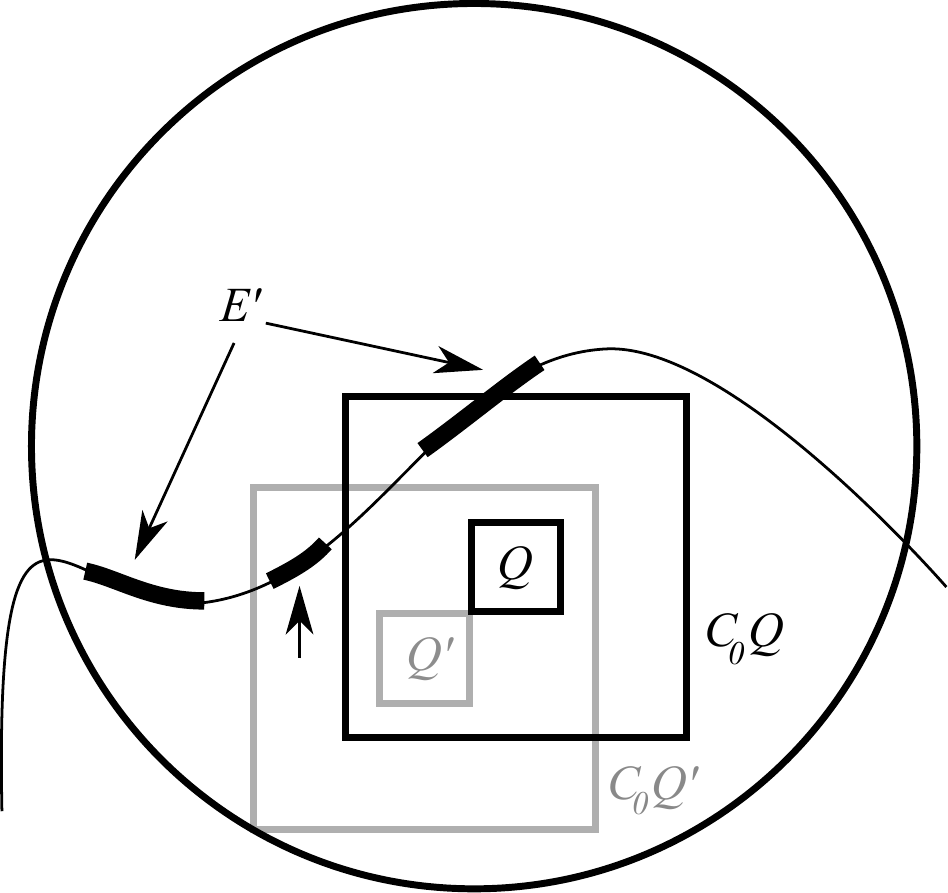}}
\caption{The cubes $Q,Q'$ and $\Delta'$ in the case that $Q\in T_{2}$.}\label{f:qq}
\end{figure}

Our goal now is to verify that $y_{Q},Q,E,E',g,$ and $Z$ satisfy the conditions of \Lemma{Zlemma}. Since $y_{Q}\in \Delta'\subseteq C_{0}Q'$ and $Q'$ is adjacent to $Q$
\[\ell(Q)\leq \dist(Q,\d\Omega)\leq \dist(y_{Q},Q)\leq \diam C_{0}Q'+\diam Q'\lec \ell(Q),\]
\[\dist(y_{Q},E')<2\ell(MB_{\Delta(Q)})\lec_{c_{0},M} \ell(Q)\]
and since $\Delta'\cap E'=\emptyset$ and $y_{Q}\in (1-t)\Delta'$,
\begin{equation}
\dist(y_{Q},E')\geq  t\ell(\Delta')\gec_{c_{0}} \ell(Q')\sim_{d} \ell(Q).
\label{e:T2}
\end{equation}

Moreover, $Q\cap B(\xi,r)\neq\emptyset$ for all $Q\in C_{E}(\tilde{\Delta})\subseteq \d\cC_{E'}^{-}(\xi,r)$, and  so by \Lemma{Zlemma},
\[\sum_{Q \in T_{2}}\ell(Q)^{d} \lec r^{d}.\]

\end{proof}

\begin{lemma}
The inequality \eqn{QindC} holds for $\d\cC_{E'}^{+}$.
\end{lemma}

\begin{proof}
The proof is basically the same as in the previous lemma, but with some minor adjustments. Let $\xi\in E'$ and $r>0$. Again, without loss of generality, it suffices to prove \eqn{QindC} for $r\leq 3r_{0}$. To see this, assume it's true and let $r>3r_{0}$, pick $n\geq 0$ so that $3^{n+1}r_{0}\geq r> 3^{n}r_{0}$. If $Q\in \cC_{j}:=\d\cC_{E'}^{+}(\xi,3^{j+1}r_{0})\backslash \d\cC_{E'}^{+}(\xi,3^{j}r_{0})$, then 
\[\ell(Q)\leq \dist(Q,E') \leq \dist(\xi,Q)\leq 3^{n+1}r_{0},\] 
yet since $E'\subseteq B(\xi_{0},r_{0})$, 
\[\ell(Q)\sim_{d,K} \dist(Q,E')\geq (3^{n}-2)r_{0}\sim 3^{n} r_{0}.\] 
Thus, all cubes in $\cC_{j}$ are contained in a ball centered about $\xi_{0}$ of radius comparable to $3^{n}r_{0}$ and have sidelengths comparable to $3^{n}r_{0}$. Hence,
\begin{multline*}
\sum_{Q\in \d\cC_{E'}^{+}(\xi,r)}\ell(Q)^{d}
\leq \sum_{Q\in \d\cC_{E'}^{+}(\xi,3^{n+1}r_{0})} \ell(Q)^{d}\\
 \leq \sum_{Q\in \d\cC_{E'}^{+}(\xi,3r_{0})}\ell(Q)^{d}+ \sum_{j=1}^{n}\sum_{Q\in \cC_{j}} \ell(Q)^{d}
 \lec_{d,K}r_{0}^{d}+ \sum_{j=0}^{n} (3^{n}r_{0})^{d}\lec 3^{nd}r_{0}^{d}\lec_{d} r^{d}.
\end{multline*}
Thus, we can assume $r\leq 3r_{0}$.

Define $\cD(\xi,r)$ just as in \Lemma{-case}. Note that if $Q\in \cC_{E}^{+}(\xi,r)$, then since $Q\in \cC_{E}^{+}\subseteq \cW_{K}(\bR^{d+1}\backslash E)$,
\[\ell(Q)\leq \dist(Q,E')\leq \dist(\xi,Q)<r\leq 3r_{0}<\ell(\Delta_{0}).\]
Moreover, we know that the parent of $Q$, $Q^{1}$, satisfies $KQ^{1}\cap E\neq\emptyset$. Thus, there is a maximal cube $\Delta(Q) \subseteq \Delta_{0}$ such that $\Delta(Q)\cap KQ^{1}\cap E'\neq\emptyset$ and $c_{0}\ell(\Delta(Q))\leq \ell(Q)$, so again $\Delta(Q)\in \cD(\xi,r)$. Again, $\# \cD(\xi,r)\lec_{C}1$, and so it suffices to show \eqn{CEDtwid}, where now
\[C_{E'}(\tilde{\Delta})=\{Q\in \d\cC_{E'}^{+}(\xi,r): \Delta(Q)\subseteq \tilde{\Delta}\}.\]

Split $C_{E'}(\tilde{\Delta})$ into sets $T_{1}$ and $T_{2}$ as before. Again, \eqn{T1} holds for the same reasons, so we're just left with estimating the sum over $T_{2}$. 

For each $Q\in C_{E'}(\tilde{\Delta})\subseteq \cC_{E'}^{+}$, we have $Q\cap \d\Omega\neq\emptyset$ by definition, so we can pick $x_{Q}\in Q\cap \d\Omega$ so that $B(x_{Q},\ell(Q))\cap \d\Omega\subseteq 3Q$. Since $\ell(Q)\sim \ell(\Delta(Q))$, we can pick $M$ large enough (depending on $d$, $K$, and $c_{0}$) so that $MB_{\Delta(Q)}\supseteq 3Q\supseteq B(x_{Q},\ell(Q))$. If $\delta$ is chosen small enough (depending on $M$ and $c_{0}$), we can guarantee that there is $y_{Q}\in E\cap B(x_{Q},\ell(Q))$. Moreover, since $MB_{\Delta(Q)}\supseteq 3Q\ni y_{Q}$, and $Q\in \cW_{K}(\bR^{d+1}\backslash E')$,
\[ \dist(y_{Q},E')\sim_{d,K} \ell(Q)\leq \dist(Q,E')\leq \diam MB_{\Delta(Q)}\lec \ell(Q),\]
and $\dist(y_{Q},Q)\leq \diam Q$ since $y_{Q}\in 3Q$. Thus, we can apply \Lemma{Zlemma} again with respect to the set $E'$.
\end{proof}

This ends the proof of \Lemma{lqsum2}.

\end{proof}

As another corollary, we get the following well known fact.

\begin{lemma}
Let $\Omega\subseteq \bR^{d+1}$ be a $C$-NTA domain with $A$-Ahlfors regular boundary. Let $E\subseteq \d\Omega\cap B(\xi_{0},r_{0})$ be a compact set with $\xi_{0}\in \d\Omega$ and $r_{0}<\diam\d\Omega$. Then there is a $C''$-NTA domain $\Omega_{E}\subseteq \Omega$ with $A''$-Ahlfors regular boundary so that $\d \Omega_{E}\cap \d \Omega=E$ and $\diam \d\Omega_{E}\geq r_{0}/C''$, where $A'',C''>0$ depend only on $A,C,$ and $d$.
\label{l:anta}
\end{lemma}

\begin{proof}
We'll just sketch some of the details. Let $\cD$ be the cubes for $\d\Omega$ and $\Omega_{E}=\Omega_{E}^{-}$ from \Lemma{omegaE-}, but pick $C$ large enough in that lemma so that, for all $Q\in \cW(\Omega)$, $C_{0}Q$ contains a cube $\Delta\in \cD$ with $c_{0}\ell(\Delta)\leq \ell(Q)< \ell(\Delta)$ and let $y_{Q}$ be the center of this cube. Following a similar procedure in the proof of \Lemma{lqsum2}, we can show that $y_{Q}$, $E$, $Z=\d \Omega$, and $\d \cC=\widetilde{\cC_{E}}^{-}(\xi,r)$ satisfy the conditions of \Lemma{Zlemma}, and so now the result follows from \Lemma{lqsum}.
\end{proof}

\section{Proof of \Theorem{uniform-main}}
\label{s:in-and-out}

\def\bbeta{b\beta}

We begin by recalling some theory from \cite{DS} and \cite{of-and-on}. 

\begin{definition}\label{d:ur}
An $A$-Ahlfors regular set $Z\subseteq \bR^{d+1}$ is {\it uniformly rectifiable} if there are constants $L,c>0$ such that, for all $\xi\in Z$ and $r\in(0,\diam Z)$, there is $E\subseteq B(\xi,r)\cap Z$ with $\cH^{d}(E)\geq cr^{d}$ and an $L$-bi-Lipschitz embedding $g:E\rightarrow \bR^{d}$.
\end{definition}

For example, if $Z$ is a bi-Lipschitz image of $\bR^{d}$, then it is trivially uniformly rectifiable. There are several different equivalent definitions of this term; for example, \cite{DS} presents seven characterisations, and in \cite{of-and-on} several more. The characterisation that will be most convenient for us, though, is one given in terms of {\it bilateral $\beta$-numbers}: for a set $Z\subseteq \bR^{d+1}$, $\xi\in Z$, $r>0$, and a hyperplane $P$ passing through $\xi$, set 
\[
\bbeta_{Z}(\xi,r,P)=\sup_{\zeta\in B(\xi,r)\cap Z}\dist(\zeta,P)/r + \sup_{\zeta\in B(\xi,r)\cap P} \dist(\zeta,Z)/r.\]
Note that by the local compactness of the Grassmanian and the continuity of $\bbeta(\xi,r,P)$ in $P$, there exists $P_{\xi,r}$ that infimizes $\bbeta(\xi,r,P)$, and we define
\[\bbeta_{Z}(\xi,r)=\bbeta_{Z}(\xi,r,P_{\xi,r}).\]

\begin{theorem} \cite[Theorem 2.4]{of-and-on}
If $Z$ is an $A$-Ahlfors regular set in $\bR^{d+1}$, then $Z$ is uniformly rectifiable if and only if, for all $\ve>0$, the set
\[\cB_{\ve}=\{(\xi,r)\in Z\times (0,\infty): \bbeta_{Z}(\xi,r)>\ve\}\]
is a {\it Carleson set}, meaning that, for all $\xi_{0}\in Z$ and $r_{0}>0$, if we define $d\sigma=d\cH^{d}|_{Z}\times \frac{dr}{r}$, then
\begin{equation}
\sigma(\cB_{\ve}\cap (B(\xi_{0},r_{0})\times (0,r_{0})))\leq C_{UR}r_{0}^{d}
\label{e:carl}
\end{equation}
where $C_{UR}$ depends on $L,d,$ and $c$ in the definition of uniform rectifiability and vice versa. 
\end{theorem}

We will say that $Z$ is {\it $C_{UR}$-uniformly rectifiable} if it satisfies \eqn{carl}.

The original definition of $\bbeta_{Z}$ infimizes over all hyperplanes $P$, not just the ones passing through $\xi$, but it's easy to see that this quantity is comparable to our current definition by a factor of two. Using these {\it centered bilateral $\beta$-numbers} will make things a bit more convenient below.

\begin{lemma}
Let $Z$ be a set in $\bR^{d+1}$, $\Sigma$ a closed set whose complement is the disjoint union of two $C$-uniform domains $\Omega_{\pm}$. Let $E=Z\cap \Sigma$, $\xi\in E$, and $r>0$. If $\bbeta_{Z}(\xi,r)<\ve<\frac{1}{8C^2}$ and $\zeta\in B(\xi,\frac{r}{2C})\cap \Sigma$ is such that $\dist(\zeta,E)>(2C+1)\ve r$, then there is $z\in Z\cap B(\xi,r)$ with $\dist(z,E)\geq \ve r$. 
\label{l:bblemma}
\end{lemma}

\begin{proof}
Let $P=P_{\xi,r}$ and $\nu$ a unit normal vector to $P$. Set 
\[H^{\pm}=\{\xi+x: \pm x\cdot \nu>0\}\] 
so $P^{c}=H^{+}\cup H^{-}$. Set $\xi_{\pm}=\xi \pm \frac{r}{4C}\nu\in H^{\pm}$. 

Let $\zeta'\in P$ be closest to $\zeta$ and let $\zeta''\in E$ be closest to $\zeta'$. Since $\bbeta_{Z}(\xi,r)<\ve$ and $\zeta'\in B(\xi,\frac{r}{2C})\cap P$
\[\dist(\zeta,P)=|\zeta-\zeta'|\geq |\zeta-\zeta''|-|\zeta''-\zeta'|\geq (2C+1)\ve r - \ve r = 2C\ve r.\]

In particular, $\zeta\not\in P$, so without loss of generality, we can suppose $\zeta\in H^{+}$. By \Definition{uniform}, we can find a curve $\gamma$ (contained in either $\Omega_{+}$ or $\Omega_{-}$) containing $\zeta$ and $\xi_{-}$ such that 
\[
\cH^{1}(\gamma)\leq C|\zeta-\xi_{-}| \leq C(|\zeta-\xi|+|\xi-\xi_{-}|)\leq C\ps{\frac{r}{2C}+\frac{r}{4C}}=\frac{3r}{4}\]
and for all $t\in \gamma$, $\dist(t,\Sigma)\geq \dist(t,\{\zeta,\xi_{-}\})/C$ (see Figure \ref{f:uniform}). 

\begin{figure}[h]
\scalebox{.57}{\includegraphics{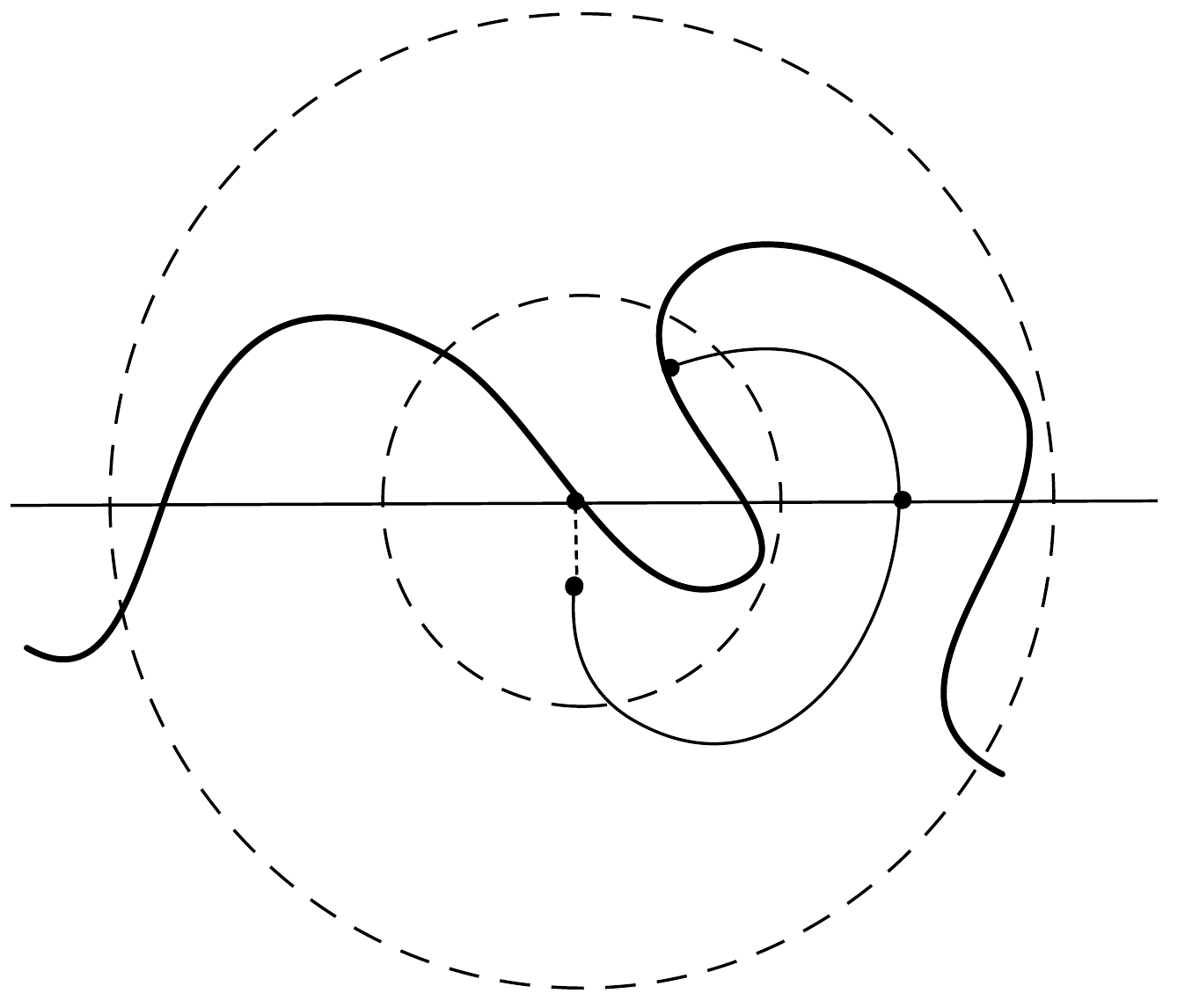}}
\begin{picture}(0,0)(300,0)
\put(180,105){$\xi$}
\put(167,80){$\xi_{-}$}
\put(210,40){$\gamma$}
\put(247,100){$t$}
\put(187,121){$\zeta$}
\put(60,90){$P$}
\put(65,65){$\Sigma$}
\put(150,140){$B(\xi,\frac{r}{2C})$}
\put(85,170){$B(\xi,r)$}
\end{picture}
\caption{Since $\zeta$ is far from $P$, we can find a point $t\in P$ far from $\Sigma$. Since $\bbeta_{Z}(\xi,r)$ is small, we can find a point in $Z$ near $t$ that will be far from $\Sigma$ as well.}\label{f:uniform}
\end{figure}
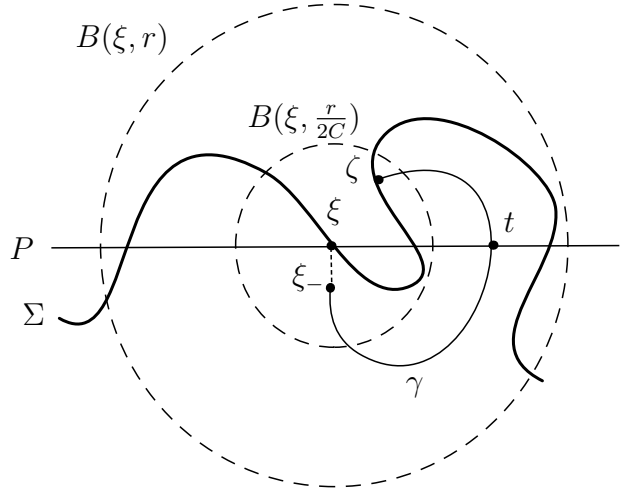

Note that $\diam\gamma\leq \cH^{1}(\gamma)\leq 3r/4$, and since $|\xi-\xi_{-}|=\frac{r}{4C}\leq r/4$, we have 
\[\gamma\subseteq B(\xi_{-},3r/4)\subseteq B(\xi,r).\]
Since $\zeta\in H^{+}$ and $\xi_{-}\in H^{-}$, there is $t\in \gamma\cap P\cap B(\xi,r)$ such that
\[\dist(t,E)\geq \dist(t,\Sigma)\geq \dist(t,\{\zeta,\xi_{-}\})/C\geq \min\{2C\ve r,\frac{r}{4C}\}/C=2\ve r\]
since $\ve<\frac{1}{8C^2}$. Since $\bbeta_{Z}(\xi,r)<\ve$, there is $z\in Z$ such that $|t-z|<\ve r$, and so
\[\dist(z,E)\geq \dist(t,E)-|t-z|\geq 2\ve r-\ve r=\ve r .\]
\end{proof}

\begin{theorem}
Let $\Omega\subseteq \bR^{d+1}$ be a $C$-uniform domain so that $(\Omega^{c})^{\circ}$ is also $C$-uniform. Let $E\subseteq \d\Omega\cap B(\xi_{0},r_{0})\cap Z$ where $Z$ is an $A$-Ahlfors regular $C_{UR}$-uniformly rectifiable set. Then there are $C^{\pm}$-NTA domains $\Omega_{E}^{\pm}$ with $A^{\pm}$-Ahlfors regular boundaries so that $\Omega_{E}^{-}\subseteq \Omega\subseteq \Omega_{E}^{+}$, $\diam \d\Omega_{E}^{\pm}\geq \diam \d\Omega /C^{\pm}$, and $E\subseteq \d\Omega_{E}^{\pm}\cap \d\Omega$. Moreover, if $\omega=\omega_{\Omega}^{z_{0}}$ where $B(z_{0},r_{0}/C)\subseteq B(\xi_{0},r_{0})\cap \Omega$, then $\omega$ is $A_{\infty}$-equivalent to $\cH^{d}$ on $E$, so
in particular, $\omega|_{E}\ll \cH^{d}|_{E} \ll \omega|_{E}$.
\label{t:in-and-out}
\end{theorem}

\begin{proof}
By \Lemma{lqsum}, it suffices to show for $\xi\in E$ and $r>0$,
\begin{equation}
\sum_{Q\in \d\cC_{E}^{\pm}(\xi,r)}\ell(Q)^{d}\lec r^{d}
\label{e:inoutCsum}
\end{equation}
where $\d\cC_{E}^{\pm}$ are defined in \eqn{dhatce-} and \eqn{dhatce+}.

We begin with $\d\cC_{E}^{-}(\xi,r)$. For $Q\in \cC_{E}^{-}$, since $C_{0}Q\cap E\neq\emptyset$, we can select $x_{Q}\in E\cap C_{0}Q$. Let $M>1$, $\ve\in (0,\min\{\frac{1}{(2C+1)M},\frac{1}{8C^{2}}\})$, and set
\[T_{1}=\{Q\in \d\cC_{E}^{-}(\xi,r): \bbeta(x_{Q},M\ell(Q))\geq \ve\}, \;\; T_{2}= \d\cC_{E}(\xi,r)\backslash T_{1}.\]

Let $Q\in T_{1}$. Note that
\begin{multline*}
\ve\leq \bbeta(x_{Q},M\ell(Q))\leq 3\bbeta(y,s) \\
\mbox{ for }(y,s)\in T_{Q}:=  (B(x_{Q},M\ell(Q))\cap Z)\times (2M\ell(Q),3M\ell(Q)).
\end{multline*}
This implies $T_{Q}\subseteq \cB_{\ve/3}$. Since $\ell(Q) \leq \dist(Q,E)\leq r$ and $x_{Q}\in Q$, we also have 
\begin{multline*}
T_{Q}\subseteq B(\xi,\diam C_{0}Q+M\ell(Q)+r)\times (0,3M\ell(Q))\\
\subseteq B(\xi,M'r)\times (0,3Mr)
\end{multline*}

where $M'=M+C_{0}\sqrt{d+1}+1$. Moreover, 
\begin{equation}
\sum_{Q\in T_{1}} \one_{T_{Q}}(x,t)\lec_{d,M} 1.
\label{e:sumTQ}
\end{equation}
To see this, observe that if $(x,t)\in T_{Q_{j}}$ for some distinct cubes $Q_{1},...,Q_{N}$, then $t\sim M\ell(Q_{j})$ for all $j$, and $\dist(x,Q_{j})\leq M\ell(Q_{j})\sim t$ for all $j$, so all $Q_{j}$ are disjoint cubes of sidelights comparable to $t/M$ contained in a ball of radius comparable to $t/M$, which implies $N\lec_{d,M}1$. Thus,
\begin{multline*}
\sum_{Q\in T_{1}}\ell(Q)^{d}
\lec_{d,M} \sum_{Q\in T_{1}} \sigma(T_{Q})
\stackrel{\eqn{sumTQ}}{\lec}_{d,M} \sigma(\cB_{\ve/3}\cap (B(\xi,M'r)\times (0,3Mr)))\\\stackrel{\eqn{carl}}{\lec}_{M,d,C_{UR}} r^{d}.
\end{multline*}

For $Q\in T_{2}$, note that since $Q\in \d\cC_{E}^{-}$, there is $Q'\in \cW(\Omega)$ such that $Q\sim Q'$ and $C_{0}Q'\cap E=\emptyset$. Again, for $C_{0}$ large enough, and since $\ell(Q)\sim_{d}\ell(Q')$, we can guarantee that there is always $z_{Q}\in \d\Omega$ such that $B(z_{Q},\ell(Q)) \subseteq  C_{0}Q'$. For $M$ large enough, $B(z_{Q},\ell(Q))\subseteq B(x_{Q},\frac{M\ell(Q)}{2C})$ (recall $x_{Q}\in C_{0}Q\cap E$), and $\dist(z_{Q},E)\geq \ell(Q)\geq \ve (2C+1)M\ell(Q)$ (since $B(z_{Q},\ell(Q))\subseteq C_{0}Q'$ and $C_{0}Q'$ does not intersect $E$). Since we also have $\ve<\frac{1}{8C^{2}}$, we can use \Lemma{bblemma} with $\Sigma=\d\Omega$ to show that there is $y_{Q}\in Z\cap B(x_{Q},M\ell(Q))$ with $\dist(y_{Q},E)\geq \ve M\ell(Q)$. Since $x_{Q}\in E$, we have $\dist(y_{Q},E)\leq M\ell(Q)$. Finally, since $Q\in \cW(\Omega)$, 
\[\ell(Q)\leq \dist(Q,\d\Omega)\leq \dist(Q,E)\leq r,\] 
and so we can apply \Lemma{Zlemma} with $\cC=T_{2}$ to show that
\[\sum_{Q\in T_{2}} \ell(Q)^{d}\lec r^{d}.\]

For $\d\cC_{E}^{+}(\xi,r)$, again, the proof is the same as above except for our choice of $x_{Q}$: For $Q\in \cC_{E}^{+} (\xi,r)$, $Q\in \cW_{K}(\bR^{d+1}\backslash E)$, and so the parent $Q^{1}$ of $Q$ satisfies $KQ^{1}\cap E\neq\emptyset$, so we can pick $x_{Q}\in KQ^{1}$. Thus, $\dist(x_{Q},Q)\leq 2K\diam Q$. We define $T_{1}$ as above, but with $M>(2K+3)\sqrt{d+1}$ so that $B(x_{Q},M\ell(Q))\supseteq 3Q$. Since $Q\cap \d\Omega\neq\emptyset$ when $Q\in \cC_{E}^{+}$, there is $z_{Q}\in Q\cap \d\Omega$, and $B(z_{Q},\ell(Q))\subseteq 3Q$. The remainder of the proof is now just like the proof we had in the case of $\d\cC_{E}^{-}(\xi,r)$, and so \eqn{inoutCsum} is proven for both cases.

The last part of the theorem now follows from \Lemma{afineq}.
\end{proof}

\begin{proof}[Proof of \Theorem{uniform-main}]

We will need the following theorem.

\begin{theorem} \cite[Theorem II]{hardsard}.
Let $D\geq d\geq 1$ and  $0<\kappa<1$ be given. 
There are constants $C'=C'(d)>0$ and $M=M(\kappa,d)$ such that if
$f:\bR^{d}\rightarrow\bR^{D}$ is a $1$-Lipschitz function, then there are sets $E_{1},...,E_{M}$ such that 
\begin{equation}
\cH_{\infty}^{d}\ps{f\ps{[0,1]^{d}\backslash\bigcup_{i=1}^{M} E_{i}}}\leq C' \kappa
\label{e:inverse-theorem-1}
\end{equation}
and such that if $E_{i}\neq\emptyset$, there is $F_{i}:\bR^{d}\rightarrow\bR^D$ which is $L_{0}$-bi-Lipschitz, $L_{0}\sim_{D}\kappa^{-1}$, so that 
\begin{equation}
F_{i}|_{E_i}= f|_{E_{i}}.
\label{e:extends}
\end{equation}
\label{t:hard-sard}
\end{theorem}

Now, let $E\subseteq \d \Omega$ be as in the statement of \Theorem{uniform-main}. Let $A\subseteq [0,r_{0}]^{d}$ and $f:A\rightarrow E$ be $L$-Lipschitz. By replacing $\Omega$ with $(Lr_{0})^{-1}\Omega$ and $f(x)$ with $\frac{f(Lr_{0}x)}{Lr_{0}}$, we may assume without loss of generality that $r_{0}=L^{-1}$ (note that this scaling does not affect the Lipschitz constant of $f$ nor the ratio $\rho=\cH^{d}_{\infty}(E)/r_{0}^{d}$). By Kirszbraun's theorem, we may extend $f$ so it is defined on all of $\bR^{d}$ and is still $L$-Lipschitz. Let $F(x)=f(x/L)$ so that $F$ is a $1$-Lipschitz map and $F([0,1]^{d})\supseteq E$. Let $\kappa= \eta \cH^{d}_{\infty}(E)/C'=\frac{\eta\rho}{L^{d}C'}$ and apply \Theorem{hard-sard} to $F$ to obtain sets $E_{1},...,E_{M}$ with $M=M(\kappa,d)$ and $L_{0}$-bi-Lipschitz functions $F_{i}:\bR^{d}\rightarrow \bR^{d+1}$ satisfying \eqn{inverse-theorem-1} and \eqn{extends}, where $L_{0}\sim_{d} \kappa^{-1}$. The sets $Z_{i}=F_{i}(\bR^{d})$ are $C_{UR}$-uniformly rectifiable sets with $C_{UR}$ depending on $d$ and $L_{0}$ (or rather, $d,\eta,L,$ and $\rho$). Let $F_{i}=Z_{i}\cap E$. By \Theorem{in-and-out}, $\omega$ is $A_{\infty}$-equivalent to $\cH^{d}$ on $F_{i}$, thus $\omega$ is also $A_{\infty}$-equivalent to $\cH^{d}$ on the finite union $E'=\bigcup_{i=1}^{M}F_{i}$. Finally
\[\cH^{d}_{\infty}(E\backslash E')
\leq \cH^{d}_{\infty}\ps{F\ps{[0,1]^{d}\backslash \bigcup_{i=1}^{M}E_{i}}}\leq C'\kappa=\eta \cH^{d}_{\infty}(E).\]

\end{proof}

\begin{proof}[Proof of \Theorem{uniformtwist}]
We simply iterate using \Theorem{uniform-main} on the set $E$ to exhaust $\cH^{d}$-almost all of $E$ with rectifiable sets upon each of which $\cH^{d}$ and $\omega$ are mutually absolutely continuous.
\end{proof}

\section{Appendix}

\subsection{The proof of \Lemma{omegaE+}}

\Lemma{omegaE+} will follow from the following lemmas. Again, we assume that $E\subseteq \d\Omega\cap B(\xi_{0},r_{0})$ where $\Omega$ is $C$-uniform (unless specified otherwise), $\xi_{0}\in \d\Omega$, and $r_{0}\in (0,\diam\d\Omega)$.

\begin{lemma}\label{l:diamo+}
If $\d \Omega$ is bounded, then $\diam \Omega_{E}^{+}\lec_{d} \diam \Omega$.
\end{lemma}

\begin{proof}\label{l:ex}
If $r\geq \diam \d\Omega$, then $\d\Omega$ is a bounded set, thus if $Q\in \cC_{E}^{+}\subseteq \cW_{K}(\bR^{d+1}\backslash E)$, $\ell(Q)\leq \dist(Q,E)\leq \diam \Omega$, and since each such $Q$ intersects $\d\Omega$, the lemma follows.
\end{proof}

\begin{lemma}
If $\Omega$ is $C$-NTA, then for $K$ large enough (depending on $\lambda,C,$ and $d$), $\Omega_{E}^{+}$ satisfies the exterior corkscrew condition. In particular, $(\Omega_{E}^{+})^{c}$ contains a ball of radius $\frac{r_{0}}{4C}$. \end{lemma}

\begin{proof}
Let $\xi\in \d\Omega_{E}^{+}$ and $r\in (0,\diam \d\Omega)$. 
\begin{enumerate}
\item Suppose $\dist(\xi,E)<r/2$. Pick $\xi'\in E$ with $|\xi-\xi'|<r/2$. Let $B=B(z,\frac{r}{2C})\subseteq \Omega\backslash B(\xi',r/2)$. If $\lambda Q\cap B(\xi',r)\cap \Sigma\neq \emptyset$, then
\begin{align*}
 \frac{K-1}{2}\ell(Q)
 & \leq \dist(Q,E)\leq \diam Q+\dist(\lambda Q,E) \\
 & \leq \sqrt{d+1}\ell(Q)+\dist(\xi',Q)\leq \sqrt{d+1}\ell(Q)+r
 \end{align*}
and so $\ell(Q)\leq r/K'$ where $K'=\frac{K-1}{2}-\sqrt{d+1}$. Thus
\[\diam \lambda Q=\lambda \sqrt{d+1}\ell(Q)\leq \frac{\lambda \sqrt{d+1}}{K'}r <\frac{r}{4C}\] 
for $K'>4\lambda C$, and thus each such $\lambda Q$ is contained in a ball of radius $\frac{r}{4C}$ centered upon $\Sigma$, and since $\dist(\frac{1}{2}B,\Sigma)\geq \frac{r}{4C}$, we have $\lambda Q\cap \frac{1}{2}B=\emptyset$. Since $\frac{1}{2}B\subseteq \Omega^{c}$ as well,  
\[\frac{1}{2} B\subseteq (\Omega_{E}^{+})^{c}\cap B(\xi',r/2)\subseteq B(\xi,r)\backslash \Omega_{E}^{+}  .\]
Clearly, $\diam(\frac{1}{2}B)\sim r$, and so we've proven the lemma in this case. Observe that, since $r_{0}<\diam \Omega$, the second part of the lemma is now proven.

\item Suppose $\dist(\xi,E')\geq r/2>0$. Then $\xi\in \d \lambda Q$ for some $Q\in \d\cC_{E}^{+}$ that intersects $\Sigma$. Again, all cubes adjacent to $Q$ have comparable diameters, so for $\lambda$ close enough to one, if $R\in \d\cC_{E}^{+}$ is the dyadic cube containing $\xi$, then $R\not\in \d \cC_{E}^{+}$, $\ell(R)\sim \ell(Q)$, and $R'=R\backslash\Omega_{E}^{+}$ is a rectangular prism with edges all of length comparable to $\ell(Q)\gec \dist(\xi,E')\geq r/2$. It is not hard to see then that $B(\xi,r)\cap R'$ contains a ball of size comparable to $r$.
\end{enumerate}
If $r\in[\diam \d\Omega,\diam \d\Omega_{E}^{+})$, then the previous lemma implies \[r<\diam\d\Omega_{E}^{+}\lec_{K,d} \diam \d\Omega.\] By the previous two cases, we know that $B(\xi,\diam\d\Omega/2)$ contains a ball of radius comparable to $\diam \d\Omega/2\sim_{K,d}r$, and thus we've proven the lemma.
\end{proof}

\begin{lemma}\label{l:uni}
$\Omega_{E}^{+}$ is uniform.
\end{lemma}

\begin{proof}

%

We will establish the lemma using Definition \ref{d:uniform}. Let $x,y\in \Omega_{E}^{+}$.

\begin{enumerate}
\item First suppose $y\in \Omega$.
\begin{enumerate}[(a)]
\item If $x\in\Omega$, then this case follows since $\Omega$ is uniform.
\item If $x\not\in \Omega$, then $x\in \lambda Q$ for some $Q\in \cC_{E}^{+}$. If $\hat{Q}$ is the union of $\lambda S$ over all $S\in \cC_{E}^{+}$ that intersect $Q$ and $y\in \hat{Q}$,  then this case follows since $\lambda Q$ is uniform.
\item Suppose  $y\not\in \hat{Q}$ but $|x-y|< \ve\ell(Q)$ for some $\ve>0$ to be chosen shortly. For $\ve>0$ small enough, this must mean that $y\in \Omega$, for otherwise $y\in \lambda S$ for some $S\in \cC_{E}^{+}$ and $\ve$ small enough implies $S\cap Q\neq\emptyset$, meaning that $y\in \hat{Q}$. Let $z\in \d\Omega_{E}^{+}$ be closest to $y$, so $|y-z|\leq |x-y|<\ve\ell(Q)$. Then $z\in S$ for some $S\in \cC_{E}^{+}$, since $\dist(z,E)\geq \dist(x,E)-|x-z|\geq \frac{K-1}{2}\ell(Q)-\ve\ell(Q)>\frac{K-1}{4}>0$ for $\ve$ small enough depending on $K$. In particular, this means $\ell(S)\gec \ell(Q)$, and so for $\ve$ small enough, depending on $\lambda$, $y\in \lambda S$, a contradiction.
\item Now suppose $y\not\in \hat{Q}$, $|x-y|\geq \ve\ell(Q)$, but assume $y\in \Omega$. Let $\xi\in \d\Omega\cap Q$, and let 
\[B(x',\frac{\lambda-2}{4C}\ell(Q))\subseteq B(\xi,\frac{\lambda-2}{4}\ell(Q))\cap \Omega\subseteq \frac{\lambda}{2}Q\cap \Omega\]
be a corkscrew ball. Since $\lambda Q$ is uniform, there is a good curve $\gamma_{1}$ between $x$ and $x'$ in $\lambda Q$ and is also a good curve in $\Omega_{E}^{+}$. Furthermore, since $x'\in \Omega$, there is a good curve $\gamma_{2}$ connecting $x'$ and $y$. Let $\gamma=\gamma_{1}\cup \gamma_{2}$. Then
\[\cH^{1}(\gamma)
\leq \cH^{1}(\gamma_{1})+\cH^{1}(\gamma_{2})
\lesssim |x-x'|+|x'-y|\lesssim \ell(Q)+|x-y|\lesssim |x-y|.
\]
Let $z\in \gamma$. If $z\in B(x',\frac{\lambda-2}{4C}\ell(Q))\cap \gamma\subseteq \frac{\lambda}{2}Q$, then 
\[ \dist(z,\d\Omega_{E}^{+})\geq \frac{\lambda-2}{4}\ell(Q)\gec |z-x|.\]
If $z\not\in B(x',\frac{\lambda-2}{4C}\ell(Q))$, then either $z\in \gamma_{1}$, in which case
\[\dist(z,\d\Omega_{E}^{+})\gec \min\{|z-x|,|z-x'|\}
\gec \min\{|z-x|,\ell(Q)\}
\gec |z-x|\]
or $z\in \gamma_{2}$, in which case 
\begin{multline*}\dist(z,\d\Omega_{E}^{+})\gec \min\{|z-y|,|z-x'|\}
\gec \min\{|z-y|,\ell(Q)\}
\\ \gec \min \{|z-y|,|z-x|\}.\end{multline*}
\end{enumerate}
\item If $y\not \in \Omega$, then $y\in \lambda R$ for some $R\in \cC_{E}|^{+}$. By the previous cases we may assume $|x-y|\geq \ve \ell(R)$, and the proof of this case is similar to the previous one, but now we also find a corkscrew ball near $y$ centered at a point $y'$ and connect paths from $x$ to $x'$, $x'$ to $y'$, and $y'$ to $y$. We omit the details. 

\end{enumerate}

\end{proof}

\begin{lemma}
If $\Omega$ is $C$-NTA, then $\diam \d\Omega_{E}^{+}\sim_{C,K,\lambda} \diam \d\Omega$.
\end{lemma}

\begin{proof}
By \Lemma{ex}, $\Omega_{E}^{+}$ has exterior corkscrews, and by \Theorem{AHMNT} and \Lemma{uni}, $\Omega_{E}^{+}$ has interior corkscrews. If $\diam \d\Omega<\infty$, we may find $B\subseteq \Omega$ so that $\diam B\sim_{C} \diam \Omega$. Then $B\subseteq \Omega_{E}^{+}$ as well. Let $B'\supseteq B$ be such that $B'\subseteq \Omega_{E}^{+}$ and there is $\xi\in \d B'\cap \d\Omega_{E}^{+}$. Since $\Omega_{E}^{+}$ has the exterior corkscrew condition, we may find $B''\subseteq  B(\xi,\diam B')\backslash \Omega_{E}^{+}$ with $\diam B''\sim_{C^{+}} \diam B' \geq \diam B\sim \diam\d\Omega$. By looking where the convex hull of $B''$ and $B'$ intersects $\d\Omega_{E}^+$, this implies $\diam d\Omega_{E}^{+}\gec \diam B'\gec \diam d\Omega$. This and \Lemma{diamo+} finish the proof. \end{proof}

\subsection{Proof of \Lemma{local-dub}}

\begin{proof}
Without loss of generality, we assume $\Delta_{0}\in \cD_{0}$. We follow the proof in \cite{Christ-T(b)}. Let $t\in (0,1)$ and 
\[E=\{\xi\in \Delta_{0}: \dist(\xi,\Sigma\backslash \Delta_{0})<t\ell(\Delta_{0})\}.\]
We can assume $E\neq\emptyset$. Let $N$ be the largest integer for which $5c_{0}^{N+1}>2t\ell(\Delta_{0})$. For $\xi\in E$, there is $\xi'\in \Sigma\backslash \Delta_{0}$ such that $|\xi-\xi'|<2t \ell(\Delta_{0})$. For every $n\geq 0$ there is $\Delta_{n}\in \cD_{n}$ such that $\xi'\in \Delta_{n}$. 

\Claim The cubes $\Delta_{n+1}$ and $\Delta_{n}$ always have distinct centers for $n=0,...,N$. 
Since $\xi'\not\in \Delta_{0}$, $\Delta_{n}\not\subseteq \Delta_{0}$, then $\Delta_{n}\subseteq \Delta_{0}^{c}$ and since $c_{1}B_{\Delta_{n+1}}\subseteq \Delta_{n+1}$ and $\xi\in \Delta_{0}$,
\begin{align*}
|\zeta_{\Delta_{n}}-\zeta_{\Delta_{n+1}}|
& \geq |\zeta_{\Delta_{n}}-\xi|-|\xi-\xi'|-|\xi'-\zeta_{\Delta_{n+1}}|\\
& \geq c_{1}\ell(\Delta_{n})-2t \ell(\Delta_{0})-\ell(\Delta_{n+1})\\
& > c_{1}\ell(\Delta_{n})-5c_{0}^{N+1}-c_{0}\ell(\Delta_{n})\\
& \geq c_{1}\ell(\Delta_{n})-2\ell(\Delta_{n+1})
=(c_{1}-2c_{0})\ell(\Delta_{n})>0.\end{align*}
This proves the claim. 

For $n\leq N$, the center of $\Delta_{n}$ is also the center of a cube $\Delta_{n+1}'\in \cD_{n+1}$, and that cube thus must be disjoint from $\Delta_{n+1}$. Moreover, it contains $c_{1} B_{\Delta_{n+1}'}=c_{0}c_{1}B_{\Delta_{n}}$, and so we have 
\begin{equation}
\Delta_{n+1}\cap c_{0}c_{1}B_{\Delta_{n}}=\emptyset.
\label{e:disjoint}
\end{equation}
Let 
\[\tilde{\cD}_{n}=\{\Delta\in \cD_{n}:\xi'\in \Delta \mbox{ for some }\xi\in E\}.\]
If $\Delta\in \tilde{\cD}_{n}$ and $n\leq N$, then 
\[|\zeta_{\Delta}-\xi|\leq |\zeta_{\Delta}-\xi'|+|\xi'-\xi|<\ell(\Delta)+2\tau\ell(\Delta_{0})<2\ell(\Delta)\]
so that
\[E\subseteq \bigcup_{\Delta\in \tilde{\cD}_{n}}2B_{\Delta} \mbox{ for all  }n<N.\]
Moreover, by \eqn{disjoint}, the family $\tilde{\cD}=\bigcup_{n=0}^{N-1}\{c_{0}c_{1}B_{\Delta}:\Delta\in \tilde{\cD}_{n}\}$ form a disjoint family of balls. Finally, all cubes in $\tilde{\cD}$ are of diameters no more than $2\ell(\Delta_{0})$ and are distance at most $2t\ell(\Delta_{0})$ from $\ell(\Delta_{0})$, so in particular they are all contained in $4B_{\Delta_{0}}$ since $t<1$. All these facts imply
\begin{align*}
\mu(E)
& \leq \frac{1}{N}\sum_{n=0}^{N-1}\sum_{\Delta\in \tilde{\cD}_{n}}\mu(2B_{\Delta})
\lec_{\mu,c_{0},c_{1}} \frac{1}{N}\sum_{n=0}^{N-1}\sum_{\Delta\in \tilde{\cD}_{n}}\mu(c_{0}c_{1}B_{\Delta})\\
& =\frac{1}{N}\mu\ps{\bigcup_{\Delta\in \tilde{\cD}} c_{0}c_{1}B_{\Delta}}
\leq \frac{\mu(4B_{\Delta_{0}})}{N} \lec_{\mu,c_{1}} \frac{\mu(\Delta_{0})}{N}.
\end{align*}
By our definition of $N$, this implies that $\mu(E)\lec_{\mu,c_{0},c_{1}}(\log\frac{1}{t})^{-1}$. In particular, this also holds if we replace $\Delta_{0}$ with any cube $\Delta\subseteq \Delta_{0}$. Thus, there is $t_{1}>0$ so that for any $\Delta\subseteq \Delta_{0}$,
\begin{equation}\label{e:mut1}
\mu(\{\xi\in \Delta: \dist(\xi,\Sigma\backslash \Delta)<t_{1}\ell(\Delta)\})<\mu(\Delta)/2.\end{equation}
If $n\in \bN$ is so that $5c_{0}^{n}<t_{1}/4\leq 5c_{0}^{n-1}$ and $\Delta_{j}\in \cD_{n}$ are such that 
\[
\Delta_{j}\cap \{\xi\in \Delta_{0}:\dist(\xi,\Sigma\backslash\Delta_{0})<t_{1}/2\}\neq\emptyset\]
then for each such $j$, since $\diam \Delta_{j}\leq \diam B_{\Delta_{j}}= 10c_{0}^{n}<t_{1}/2$,
\begin{equation}\label{e:djinring}
\Delta_{j}\subseteq  \{\xi\in \Delta_{0}:\dist(\xi,\Sigma\backslash \Delta_{0})<t_{1}\}.\end{equation}
Suppose we have shown for some $m\geq 1$ that for any $\Delta\subseteq \Delta_{0}$
\begin{equation}\label{e:mutm}
\mu (\{\xi\in \Delta: \dist(\xi,\Sigma\backslash \Delta)<c_{0}^{n(m-1)}t_{1}\ell(\Delta)/2\})
<2^{-m}.\end{equation}
(Note that the $m=1$ case follows from \eqn{mut1}.) Then, recalling that $\Delta_{0}\in \cD_{0}$,
\begin{align*}
\mu (\{\xi\in &\Delta_{0}  : \dist(\xi,\Sigma\backslash \Delta_{0})<c_{0}^{mn}t_{1}\ell(\Delta_{0})/2\})\\
& \leq \sum_{j} \mu(\{\xi\in \Delta_{j}: \dist(\xi,\Sigma\backslash\Delta_{j})<c_{0}^{n(m-1)}t_{1} \ell(\Delta_{j})/2\})\\
& \stackrel{\eqn{mutm}}{<}2^{-m} \sum_{j}\mu (\Delta_{j})
\stackrel{\eqn{djinring}}{\leq}2^{-m}\mu(\{\xi\in \Delta_{0}: \dist(\xi,\Sigma\backslash \Delta_{0})<t_{1}\ell(\Delta_{0})\}\\
& \stackrel{\eqn{mut1}}{<}2^{-m-1}\mu(\Delta_{0}).
\end{align*}
By induction, \eqn{mutm} holds for all $m\geq 1$, which finishes the proof.

\end{proof}

\bibliographystyle{amsplain}

\def\cprime{$'$}
\providecommand{\bysame}{\leavevmode\hbox to3em{\hrulefill}\thinspace}
\providecommand{\MR}{\relax\ifhmode\unskip\space\fi MR }
\providecommand{\MRhref}[2]{%
  \href{http://www.ams.org/mathscinet-getitem?mr=#1}{#2}
}
\providecommand{\href}[2]{#2}

\end{document}